\documentclass[aos,MSNbibl,nameyear,dvips]{arximspdf}
\usepackage{accents}
\usepackage{graphicx}
\doi{10.1214/12-AOS1053} 
\volume{40}
\issue{6}
\pubyear{2012}
\firstpage{2850}
\lastpage{2876}

\makeatletter
\renewcommand{\mathring}[1]{\accentset{\circ}{#1}}
\renewcommand{\P}{\mathrm{P}}
\newtheorem{theorem}{Theorem}
\newtheorem{lemma}{Lemma}
\newtheorem{corollary}{Corollary}

\newproclaim{remark}{Remark}
\makeatother

\begin{document}
\begin{frontmatter}

\title{Optimal two-stage procedures for estimating location and size
of the maximum of a~multivariate regression function}
\runtitle{Optimal two-stage procedure}

\begin{aug}
\author[A]{\fnms{Eduard} \snm{Belitser}\corref{}\ead[label=e1]{e.n.belitser@tue.nl}},
\author[B]{\fnms{Subhashis} \snm{Ghosal}\thanksref{m2}\ead[label=e2]{sghosal@stat.ncsu.edu}}
\and
\author[C]{\fnms{Harry} \snm{van Zanten}\ead[label=e3]{j.h.vanzanten@uva.nl}}
\thankstext{m2}{Supported by a
grant from NWO of the Netherlands.
Part of this work was done when
S. Ghosal was visiting EURANDOM, Eindhoven.}
\runauthor{E. Belitser, S. Ghosal and H. van Zanten}
\affiliation{Eindhoven University of Technology,
North Carolina State University and~University~of~Amsterdam}
\address[A]{E. Belitser\\
Department of Mathematics\\
Eindhoven University of Technology\\
P.O. Box 513\\
5600 MB Eindhoven\\
The Netherlands\\
\printead{e1}}

\address[B]{S. Ghosal\\
Department of Statistics\\
North Carolina State University\\
4276 SAS Hall, 2311 Stinson Drive\\
Raleigh, North Carolina 27695-8203\\
USA\\
\printead{e2}}

\address[C]{H. van Zanten\\
Korteweg-de Vries Institute for Mathematics\\
University of Amsterdam\\
P.O. Box 94248\\
1090 GE Amsterdam\\
The Netherlands\\
\printead{e3}}
\end{aug}

\received{\smonth{9} \syear{2011}}
\revised{\smonth{9} \syear{2012}}

\begin{abstract}
We propose a two-stage procedure for estimating the location
$\bolds{\mu}$ and size~$M$ of the maximum of a smooth $d$-variate
regression function $f$. In the first stage, a preliminary estimator
of $\bolds{\mu}$ obtained from a standard nonparametric smoothing method
is used. At the second stage, we ``zoom-in'' near the vicinity of the
preliminary estimator and make further observations at some design
points in that vicinity. We fit an appropriate polynomial regression
model to estimate the location and size of the maximum.
We establish that, under suitable smoothness conditions and appropriate choice
of the zooming, the second stage estimators have better convergence
rates than the corresponding first stage estimators of $\bolds{\mu}$ and
$M$.  More specifically, for $\alpha$-smooth regression functions, the
optimal nonparametric rates $n^{-(\alpha-1)/(2\alpha+d)}$ and
$n^{-\alpha/(2\alpha+d)}$ at the first stage can be improved to
$n^{-(\alpha-1)/(2\alpha)}$ and $n^{-1/2}$, respectively, for
$\alpha>1+\sqrt{1+d/2}$. These rates are optimal in the class of all
possible sequential estimators. Interestingly, the two-stage procedure
resolves ``the curse of the dimensionality'' problem to some extent, as the
dimension $d$ does not control the second stage convergence rates,
provided that the function class is sufficiently smooth.
We consider a multi-stage generalization of our procedure that
attains the optimal rate for any smoothness level $\alpha>2$
starting with a preliminary estimator with any power-law rate at the first
stage.
\end{abstract}



\end{frontmatter}

\section{Introduction}\label{sec1}

In many applications, it is of interest to estimate the location
and size of the extremum of a univariate or multivariate regression
function. For instance,\vadjust{\goodbreak} an oil company may
be interested in determining  the best location for drilling a well in a
confined region. Based on information obtained from drilling
at a few preliminary locations in the region, the goal is to obtain an
estimate of the best location and the amount of the reserve based on
these noisy measurements.

Suppose we observe noisy measurements of an unknown regression
function
$f\dvtx \mathbb{R}^d \to \mathbb{R}$, sampled at
points from some compact, convex set $D\subset \mathbb{R}^d$,
\begin{equation}
\label{data} Y_k= f(\mathbf{x}_k)+\xi_k,\qquad
\mathbf{x}_k \in D \subset \mathbb{R}^d, \qquad k=1,\ldots,n,
\end{equation}
where the $\xi_k$'s are independent zero mean errors with
$\operatorname{Var}(\xi_k) = \sigma^2$. Clearly, for estimating any feature of
$f$, the estimation error increases with $\sigma^2$. Thus among all
error distributions satisfying $\operatorname{Var}(\xi_k)\le \sigma^2$ for
every $k$, the homoscedasticity condition $\operatorname{Var}(\xi_k) =
\sigma^2$ is the least favorable. This shows that the latter condition
can be relaxed to the former without increasing the bound on error of
estimation, and the obtained rates under the homoscedasticity
condition remains minimax optimal under the larger
heteroscedastic model.

Assume that $f$ has a unique maximum  at $\bolds{\mu}$ in
the interior of $D$, that is,
\begin{equation}
\label{new1} \max_{\mathbf{x}\in D} f(\mathbf{x}) = f(\bolds{\mu})=M,\qquad f(
\mathbf{x})<f(\bolds{\mu})\qquad \mbox{for all }\mathbf{x}\ne \bolds{\mu}.
\end{equation}
If the function $f$ is sufficiently smooth, then
the gradient $\nabla f(\bolds{\mu})=0$ and the Hessian matrix
of $f$ at $\bolds{\mu}$ is nonpositive definite. The goal is to estimate
the maximum of the regression function
$M=f(\bolds{\mu})$ and its location $\bolds{\mu}$.

Clearly, the choice of the design points  $\{\mathbf{x}_k$, $k=1,\ldots,n\}$
significantly influences the estimation accuracy. There are two basic
design settings:
fixed in advance (or randomly sampled from a chosen distribution)
and sequential, where one is allowed to use the information obtained
from an earlier sample to determine subsequent design points. If the
design is fixed and nothing is known
about the location of the maximum,  the design points should be
``almost uniformly'' spread out all over the set of interest $D$.
The problem of estimating the location and size of extrema
of nonparametric regression functions for the fixed design situation
has been studied by many authors.
The one-dimensional case is thoroughly investigated, whereas the study in
the multivariate situation has been limited; see
M\"uller (\citeyear{Muller1985,Muller1989}), \citet{ShoungZhang2001},
Facer and M\"uller (\citeyear{FacerMuller2003}) and the references therein.
The minimax rate for estimating the maximum value of the function ranging
over an $\alpha$-smooth nonparametric class (e.g.,
isotropic H\"older class defined below) is $n^{-\alpha/(2\alpha+d)}$.
As to the estimation of the location of the maximum, it is a folklore that
the minimax rate is the same as the minimax rate for estimating the
first derivative of the regression function, which is given by
$n^{-(\alpha-1)/(2\alpha+d)}$. In the setting of estimating the mode
$\mu$ of a univariate twice differentiable density $f$, Hasminskii
(\citeyear{Hasminskii1979}) showed that under the assumption that $f''(\mu)<0$,
the lower bound for the minimax risk rate  is of the order $n^{-1/5}$
consistent with the rate $n^{-(\alpha-1)/(2\alpha+d)}$.
Klemel\"a (\citeyear{Klemela2005}) considered the problem of adaptive estimation of the mode
of a multivariate density with a bounded support that satisfies, in
a neighborhood of the mode, a smoothness condition of a level higher
than $2$.

If we can choose a design point before making each observation using
the data obtained so far, then we are in the classical sequential
design setting. Kiefer and Wolfowits (\citeyear{KieferWolfowitz1952}) introduced
a Robbins--Monro type of algorithm to estimate the mode $\mu$ of $f$ in
the univariate framework. \citet{Blum1954} proposed a multivariate version
of their algorithm which allows to estimate
the location $\bolds{\mu}$ of the maximum of a multivariate regression function $f$.
Since then, this Kiefer--Wolfowits--Blum recursive algorithm
has been extended in many directions by many authors.
The main fact is that the algorithm converges to $\bolds{\mu}$ with the
rate $n^{-1/3}$ under the assumption that the regression function $f$
is three times differentiable.
More generally, \citet{Chen1988} and \citet{PolyakTsybakov1990}
established that, in the  sequential design setting,
the minimax rate for estimating the location of the maximum
of $\alpha$-smooth regression functions is $n^{-(\alpha-1)/(2\alpha)}$.
\citet{Dippon2003} proposed a general class of randomized gradient
recursive algorithms which attain the optimal convergence rate.
\citet{MokkademPelletier2007} considered the problem of simultaneously
estimating, in the  sequential design setting, the location and the
size of the maximum of a regression function that is three times
continuously differentiable.
They proposed a companion recursive procedure to
the Kiefer--Wolfowits--Blum algorithm so that, by applying
both the companion and the Kiefer--Wolfowits--Blum algorithms,
one can simultaneously estimate the location and size of the
maximum of regression functions in an on-line regime.
Interestingly, in a sequential design setting, the convergence
rate for estimating the maximum itself $M=f(\bolds{\mu})$
can, in principle, attain the parametric rate $n^{-1/2}$.
The companion procedure of \citet{MokkademPelletier2007}
for estimating the maximum can also achieve the parametric
rate $n^{-1/2}$, but this companion procedure must use different design points
than those used in the Kiefer--Wolfowits--Blum procedure.

In this paper, we propose a two-stage strategy to tackle the
problem of simultaneously estimating the location $\bolds{\mu}$ and size
$M$ of the maximum of the regression function $f$ according to the
observation scheme (\ref{data}). This is an approach in between the two
above described  frameworks---global fixed design and a
fully sequential design. Often, from an operational point of view,
fully sequential sampling can be expensive, whereas a two-stage
procedure is much simpler to implement.
Our findings establish that the
two-stage procedure can be properly designed to match the strength of
a fully sequential procedure. Moreover, the same design scheme can be
used to obtain the optimal rates for estimating both $\bolds{\mu}$ and
$M$.

Now we describe the two-stage procedure. We construct a preliminary
estimator $\tilde{\bolds{\mu}}$ of $\bolds{\mu}$ by spending a portion of
our sampling budget to make observations over a relatively uniform
grid of points in the area of interest and applying some standard
nonparametric smoothing method for the fixed design setting based on
this initial set of data. Additional prior information, if available,
may also be used to reduce the span of the design points or to more
efficiently choose design points leading to increased accuracy of the
preliminary estimator. At the second stage, we ``zoom-in'' on a
neighborhood of $\tilde{\bolds{\mu}}$ of an appropriate size $\delta_n$,
to be called the \textit{localization parameter}. The idea is that if
this vicinity is ``small enough,'' that is, the preliminary estimator
$\tilde{\bolds{\mu}}$ converges to $\bolds{\mu}$, the regression function
$f$ can be accurately approximated by a Taylor polynomial within the
vicinity of $\tilde{\bolds{\mu}}$. We then spend the remaining portion
of the sampling budget to gather further observations at appropriately chosen
design points in the vicinity of $\tilde{\bolds{\mu}}$.
Finally, we fit a polynomial regression model on the new set of data
and show that the remainder of the expansion is appropriately small,
provided that the preliminary estimator $\tilde{\bolds{\mu}}$ has
sufficient accuracy. This procedure leads to improved estimators of
$\bolds{\mu}$ and $M$ and does not use knowledge of the noise variance
$\sigma^2$. The last step in our approach is reminiscent of the
nonparametric methodology of local polynomial regression in case of
fixed design setting; see \citet{FanGijbels1996}. Our two-stage
procedure is motivated by the recent work of \citet{LanBanerjeeMichailidis2011} and
\citet{TangBanerjeeMichailidis2011}, who, respectively, considered such procedures for estimating change
points in a regression function and the level point of a univariate monotone
regression function. Motivating grounds for a two-stage approach were
nicely described by them. The principal differences between their and
our techniques are that we consider smooth rather than step or monotone
functions, and we use polynomial regression of an appropriate
degree in the second stage rather than regression based on step or
linear functions, respectively, used by them.

The results for estimating $\bolds{\mu}$ and $M$ under the fully
sequential setting, which we are aware of, all follow the
Robbins--Monro procedure, where the next design point depends only on
the previous observation and does not incorporate all available
information up to the current moment. In this setting, one makes
observations only along a certain path of design points, eventually
leading to the location of the maximum. In our two-stage approach, one
also gets the global estimate of the regression function from the
first stage all over the area of interest, which may be useful in some
practical situations. We also get an accompanying estimator for the
size of the maximum $M$ (in fact, for all the relevant derivatives at
the location of the maximum) in a natural way, while in a
Robbins--Monro type sequential design, one needs to adjust the design
points to estimate $M$. This can place serious constraints on the
available budget since typically both $\bolds{\mu}$ and $M$ need to be
estimated.

Our main result gives a decomposition of the convergence
rate of the second stage estimator as the sum of an approximation
term and a stochastic  term,
similar to the classical bias-variance trade-off.
An implication of the main result is as follows.
Suppose we take a preliminary
nonparametric estimator $\tilde{\bolds{\mu}}$  with the optimal
single-stage convergence rate $n^{-(\alpha-1)/(2\alpha+d)}$.
Then  by applying  our two-stage procedure with an appropriate
choice of the localization parameter $\delta_n$,
we obtain optimal (for the sequential design setting)
convergence rates, $n^{-(\alpha-1)/(2\alpha)}$ and $n^{-1/2}$, respectively,
under the condition on the smoothness parameter $\alpha >1+\sqrt{1+d/2}$.
Note that $n^{-1/2}$ is also the ``oracle rate'' for estimating
$M$ corresponding to taking $n$ samples at the ``perfect location''
$\bolds{\mu}$.
Thus, for $\alpha$-smooth regression functions, the second stage
improves the rates in estimating $\bolds{\mu}$ and $M$
from the nonparametric rates $n^{-(\alpha-1)/(2\alpha+d)}$ and
$n^{-\alpha/(2\alpha+d)}$ to the optimal sequential rates
$n^{-(\alpha-1)/(2\alpha)}$ and $n^{-1/2}$, respectively.
Curiously, the dimension $d$ disappears from powers in the second stage
convergence rates. However, the curse of dimensionality is still
present in a milder form through the constraint
$\alpha >1+\sqrt{1+d/2}$. For instance, if $\alpha>3$, then the second
stage rates are optimal for $d=1,\ldots, 6$.
We can resolve the curse of dimensionality completely by considering
a multi-stage generalization of the two-stage procedure, obtained by
iterating the second stage operation on the estimator obtained in the
second stage, and continuing the iteration sufficiently many times.
We shall show that after an appropriate number of stages, the optimal
convergence rates are attained for any $\alpha>2$.
In fact, even if we start with a not necessarily optimal preliminary
estimator at the first stage (as long as it has a convergence rate  of
a power-law type), this multi-stage approach will lead to the optimal
resulting stage after a finite number of stages. The number of stages
depends on the smoothness of the regression function
and the quality (convergence rate) of the preliminary estimator
from the first stage. The method still uses knowledge of the
smoothness level $\alpha$ in its formulation, and hence is not
adaptive for estimating $\bolds{\mu}$. Nevertheless, the multi-stage
procedure achieves the optimal rate $n^{-1/2}$ for
estimating $M$ without using the knowledge of $\alpha$.

The paper is organized as follows. In Section~\ref{sec2},
we introduce the notation and assumptions. Section~\ref{sec3} describes the
two-stage procedure and states the main result. The multi-stage
generalization is discussed in Section~\ref{sec4}, and some simulation
results are given in Section~\ref{sec5}. Proofs are presented in
Section~\ref{sec6}. Some auxiliary results are given in the \hyperref[app]{Appendix}.

\section{Notation, preliminaries and assumptions}\label{sec2}

We describe the notation and conventions to be used in this paper.
All asymptotic relations and symbols
[like $O(\delta_n)$, $o(\delta_n)$, $O_p(\delta_n)$,
$o_p(\delta_n)$  etc.] will refer to the asymptotic regime
$n\to \infty$; here $c_n=O(\delta_n)$ [resp.,
$c_n=o(\delta_n)$] means that that $c_n/\delta_n$ is bounded
(resp., $c_n/\delta_n\to 0$) and for a stochastic sequence
$X_n$, $X_n=O_p(\delta_n)$ [resp., $X_n=o_p(\delta_n)$] means
that that $\P\{|X_n|\le Kc_n\}\to 1$ for some constant $K$
(resp., $\P\{|X_n|<\varepsilon\delta_n\}\to 1$ for all $\varepsilon>0$).
For numerical sequences $\beta_n$ and $\beta'_n$,
by $\beta_n \ll \beta'_n$ (or $\beta'_n \gg \beta_n$)
we mean that $\beta_n=o(\beta'_n)$, while by $\beta'_n \gtrsim \beta_n$,
we mean that $\beta_n=O(\beta'_n)$.
By $\beta_n \asymp \beta'_n$ we mean that $\beta_n= O(\beta'_n)$
and $\beta'_n=O(\beta_n)$.
Let $\mathbb{N}$ stand for $\{0,1,2,\ldots\}$.
For a set~$S$, denote by $|S|$ the number of elements in $S$.
Vectors are represented by bold symbols and can be upper or lowercase
English or Greek letters. All vectors are in the column format with
the corresponding nonbold letters with subscripts denoting the
components, that is, for $\mathbf{x}, \mathbf{x}_k \in \mathbb{R}^d$,
$\mathbf{x}=(x_1,\ldots, x_d)$ and $\mathbf{x}_k=(x_{k,1},\ldots, x_{k,d})$.
By $\|\mathbf{x}\|$ for a vector $\mathbf{x}$, we mean the usual Euclidean
norm of $\mathbf{x}\in\mathbb{R}^d$. Matrices are also written in bold and
only uppercase English letters are used to denote them.
If~$\mathbf{A}$ is a matrix, $\|\mathbf{A}\|$ will stand for a norm on the
space of matrices such as the operator norm defined by $\|\mathbf{A}\| =
\sup_{\|\mathbf{x}\| \le 1} \|\mathbf{A}\mathbf{x}\|$.
Let $B(\mathbf{c}, R)=\{\mathbf{z}\in\mathbb{R}^d\dvtx \|\mathbf{z}-\mathbf{c}\| \le R\}$ denote a ball in $\mathbb{R}^d$ with center
$\mathbf{c}\in\mathbb{R}^d$ and radius $R>0$.
Define a cube around a point $\mathbf{a}=(a_1,\ldots, a_d) \in \mathbb{R}^d$
with an edge length $2\delta$ by
\begin{equation}
\label{multcube} C(\mathbf{a},\delta)=\bigl\{\mathbf{x} \in
\mathbb{R}^d\dvtx x_k \in [a_k- \delta,
a_k+\delta], k=1,\ldots, d\bigr\} \subset \mathbb{R}^d.
\end{equation}
If $\mathbf{a} = \mathbf{0}$, then we write $C(\delta)$ for $C(\mathbf{0},\delta)$.

We shall use the multi-index notation
$\mathbf{i} = (i_1,\ldots, i_d) \in \mathbb{N}^d$.
For a multi-index $\mathbf{i}$, a vector $\mathbf{x}\in\mathbb{R}^d$
and a sufficiently smooth function $f$ of $d$ variables, define
\[
|\mathbf{i}| = \sum_{k=1}^d
i_k,\qquad  \mathbf{i}! = \prod_{k=1}^d
i_k!, \qquad \mathbf{x}^{\mathbf{i}}= \prod_{k=1}^d
x_k^{i_k}, \qquad \mathbf{D}^{\mathbf{i}} f(\mathbf{x}_0)
= \frac{\partial^{|\mathbf{i}|} f(\mathbf{x})}{\partial x_1^{i_1}
\cdots \partial x_d^{i_d}}\Big |_{\mathbf{x}=\mathbf{x}_0}.
\]

For $k,d,r \in \mathbb{N}$, define
\[
\mathbb{I}_k(d)= \bigl\{\mathbf{i} \in \mathbb{N}^d\dvtx i_1+\cdots+i_d = k \bigr\},\qquad \mathbb{I}(r,d)=\bigcup
_{k=0}^r\mathbb{I}_k(d),
\]
with $\mathbb{I}_0(d)=\{(0,\ldots,0)\}$.
For convenience in writing, $\mathbb{I}(r,d)$ will be enumerated by
stacking elements of $\mathbb{I}_0(d), \mathbb{I}_1(d),\ldots,
\mathbb{I}_r(d)$, in that order.
Within each $\mathbb{I}_k(d)$, the elements are arranged
following the lexicographic (or dictionary) ordering.
Observe that $\mathbb{I}_k(d)$ and $\mathbb{I}_l(d)$ introduced above
are disjoint if $k\not= l$. The cardinality $|\mathbb{I}_k(d)|$ is
the number of $d$-tuples $(k_1,\ldots, k_d)\in \mathbb{N}^d$
such that $k_1+\cdots +k_d = k$, or equivalently, the number
of ways to put $k$ balls in $d$ boxes. Thus
$|\mathbb{I}_k(d)|=\bigl({d+k-1 \atop k}\bigr)=\bigl({d+k-1 \atop d-1}\bigr)$,
and hence
\[
\bigl|\mathbb{I}(r,d)\bigr|= \sum_{k=0}^r\bigl|
\mathbb{I}_k(d)\bigr|= \sum_{k=0}^{r}
\pmatrix{d+k-1
\cr
d-1}.
\]
In particular, $|\mathbb{I}(r,1)|=r+1$.

For an $\alpha\in\mathbb{R}$, let $\lceil \alpha \rceil$
be the smallest integer bigger than or equal to $\alpha$. Then
$r_\alpha=\lceil \alpha-1 \rceil$ stands for the largest integer which
is strictly less than $\alpha$. Clearly, if $\alpha\in \mathbb{N}$,
then $r_\alpha = \alpha-1$.

For $\alpha,L>0$ and a compact, convex set $D\subseteq \mathbb{R}^d$,
introduce an isotropic H\"older functional class
$\mathcal{H}_d(\alpha,L,D)$, consisting of $r_\alpha$-times
differentiable functions $f\dvtx D \to \mathbb{R}$ such that
\begin{equation}
\label{multHolder} \bigl|f(\mathbf{x})- P_{f,\mathbf{x}_0}(\mathbf{x})\bigr|
\le L \|\mathbf{x} -\mathbf{x}_0\|^{\alpha},\qquad \mathbf{x},
\mathbf{x}_0 \in D,
\end{equation}
where
\begin{equation}
\label{new3} P_{f,\mathbf{x}_0}(\mathbf{x}) = \sum_{\mathbf{i} \in \mathbb{I}(r_\alpha,d)}
\frac{1}{\mathbf{i}!} \mathbf{D}^{\mathbf{i}}f(\mathbf{x}_0) (
\mathbf{x}-\mathbf{x}_0)^{\mathbf{i}}
\end{equation}
is the Taylor polynomial of order $r_\alpha$ obtained by expansion
of $f$ about the point~$\mathbf{x}_0$.

Put $q(\alpha,d)=|\mathbb{I}(r_\alpha,d)|-1$. Observe that the total
number of terms in the $d$-variate Taylor polynomial
$P_{f,\mathbf{x}_0}(\mathbf{x})$ of order $r_\alpha$ defined in
(\ref{new3}) is $q(\alpha,d)+1$.

For a function $g\dvtx \mathbb{R}^d \to \mathbb{R}$ such that all
second-order partial derivatives of $g$ exist at a point
$\mathbf{x}_0 \in \mathbb{R}^d$, denote by $Hg(\mathbf{x}_0)$ the Hessian
matrix of the function $g$ at the point $\mathbf{x}_0$, whose $(i,j)$th
entry is given by $ \frac{\partial^2 g(\mathbf{x}_0)}{
\partial x_j\, \partial x_i}$, $i,j=1,\ldots,d$.
Notice that if $g$ has continuous second order partial derivatives
at $\mathbf{x}_0$, then the Hessian matrix $Hg(\mathbf{x}_0)$ is symmetric, and
hence its eigenvalues must be real.
For a symmetric matrix $\mathbf{M}$, denote by $\lambda_{\min}(\mathbf{M})$
and $\lambda_{\max}(\mathbf{M})$ the smallest and the largest eigenvalues
of $\mathbf{M}$, respectively.

Consider the model (\ref{data}) with $f\dvtx D\to \mathbb{R}$.
We now describe the assumptions on~$f$ to be used throughout the paper.

\begin{longlist}[(A1)]
\item[(A1)]
The function $f(\mathbf{x})$, $\mathbf{x} \in D \subseteq \mathbb{R}^d$,
allows extension on a slightly bigger set $D^\varepsilon=\bigcup_{\mathbf{x}\in
D} B(\mathbf{x},\varepsilon)$ for some $\varepsilon > 0$ (in order to avoid
boundary effects) and belongs to an isotropic H\"older functional
class $\mathcal{H}_d(\alpha,L,D^\varepsilon)$ defined by
(\ref{multHolder}), with $L>0$ and $\alpha>2$.

\item[(A2)] There is a unique point $\bolds{\mu}$ in the interior
$\mathring{D}$ of $D$ that maximizes  the function $f$ on $D$, that is,
$M=\sup_{\mathbf{x}\in D} f(\mathbf{x})=\max_{\mathbf{x}\in \mathring{D}}
f(\mathbf{x}) = f(\bolds{\mu})$ and  $\lambda_{\max}(Hf(\bolds{\mu}))<0$.
\end{longlist}

Note that conditions (A1) and (A2) imply that  $\nabla f(\bolds{\mu})=\mathbf{0}$,
and the Hessian $Hf(\bolds{\mu})$ is a symmetric
and negative definite matrix.
Besides, as $\alpha>2$, the Hessian matrix $Hf(\mathbf{x})$
is continuous and therefore for some $\kappa,\lambda_0>0$,
\begin{equation}
\label{newA3} \sup_{\mathbf{x}\in B(\bolds{\mu}, \kappa)}\lambda_{\max}\bigl(Hf(
\mathbf{x})\bigr)\le -\lambda_0.
\end{equation}
Notice that constants $\kappa$, $\lambda_0$ depend on $f$.
If we do not pursue any uniformity over~$f$ in our results, then the condition (\ref{newA3}) follows
from (A1), (A2) and can therefore be used in the
proofs. However, when uniformizing the results over a functional class,
this condition becomes autonomous and must be
added to the description of the functional class; see Remark
\ref{remark1} below.

\section{The two-stage procedure}\label{sec3}

For a column vector
${\bolds{\vartheta}}=(\vartheta_{\mathbf{i}}\dvtx \mathbf{i} \in \mathbb{I}(r_\alpha,d))^T$,
introduce the multivariate polynomial function
\begin{equation}
\label{multpolynomial} f_{\bolds{\vartheta}}(\mathbf{x}) = \sum
_{\mathbf{i}\in\mathbb{I}(r_\alpha,d)} \vartheta_{\mathbf{i}} \mathbf{x}^{\mathbf{i}}
=\sum_{k=0}^{r_\alpha} \sum
_{\mathbf{i}\in\mathbb{I}_k(d)} \vartheta_{\mathbf{i}} \mathbf{x}^{\mathbf{i}}.
\end{equation}

We now describe the two-stage procedure for estimating  the parameters
$(\bolds{\mu},M)$. The first stage consists of the first two steps and
the steps 3--5 comprise the second stage.
\begin{longlist}[(5)]
\item[(1)]
The first stage starts as follows. For $\upsilon \in
(0,1)$, choose first stage design budget, that is, $n_1\in \mathbb{N}$ such
that $0<n_1 <n$, ${n_1}/{n} \to \upsilon$.  Find $n_1$
design points $\{\tilde{\mathbf{x}}_i,   i=1,\ldots,n_1\}$
approximately uniformly over the set $D$ in the sense that,
for some $c_1,c_2 >0$, the family of balls
$\{B(\tilde{\mathbf{x}}_i, c_1 n^{-1/d}),   i=1,\ldots,n_1\}$ covers $D$ and
$\|\tilde{\mathbf{x}}_i-\tilde{\mathbf{x}}_j\| \ge  c_2 n^{-1/d}$ for $i\not=j$.

Observe the data
$\mathcal{D}_1^*=\{(\tilde{\mathbf{x}}_k,\tilde{Y}_k),   k=1,\ldots,  n_1\}$,
$\tilde{Y}_k =f(\tilde{\mathbf{x}}_k) +\tilde{\xi}_k$, $k=1,\ldots,n_1$,
according to the model
(\ref{data}).

\item[(2)]
Using $\mathcal{D}_1^*$,  construct a preliminary consistent estimator
$\tilde{\bolds{\mu}}$ of $\bolds{\mu}$. For $d=1$, one may use the kernel
estimator of M\"uller (\citeyear{Muller1989}) and for $d\ge 2$,
its multivariate generalization given by Facer and M\"uller (\citeyear{FacerMuller2003}).

\item[(3)]
Let $n_2=n-n_1$ be the remaining portion of the design budget,
and let $l$ be the smallest integer that satisfies $2l\ge r_\alpha$.
Assume that $n_2 = n_3 (2l+1)^d$ for some $n_3\in\mathbb{N}$,
which is always possible to arrange.
Note that $n_3 \ge c n$ for some constant $c>0$.
Introduce a \textit{localization parameter}
$\delta_n>0$, $\delta_n\to 0$, and define the set
\[
\prod_{i=1}^d \{\tilde{
\mu}_k+j_i\delta_n, j_i = 0,\pm
1, \pm 2, \ldots, \pm l\}= \{\tilde{\mathbf{d}}_1, \ldots, \tilde{
\mathbf{d}}_{(2l+1)^d} \},
\]
which consists of $(2l+1)^d$ different points
from  the $d$-dimensional cube $C(\tilde{\bolds{\mu}},l\delta_n)$.

Now introduce the second stage design
points $\{\mathbf{x}_k,   k=1,\ldots,n_2\}$
in such a way that $|I_j|=n_3$ for all $j=1,\ldots, (2l+1)^d$, where
$I_j = \{1\le k \le n_2\dvtx \mathbf{x}_k =\tilde{\mathbf{d}}_j\}$.
In other words, each point among the $(2l+1)^d$ different points from the
set $ \{\tilde{\mathbf{d}}_1, \ldots,
\tilde{\mathbf{d}}_{(2l+1)^d} \}$
is repeated $n_3=n_2/(2l+1)^d$ times in the second stage design
$\{\mathbf{x}_k,   k=1,\ldots,n_2\}$.
Observe the data $\mathcal{D}_2^*= \{(\mathbf{x}_k,Y_k),   k=1,\dots, n_2\}$,
$Y_k =f(\mathbf{x}_k) +\xi_k,   k=1,\ldots,n_2$, according to the model
(\ref{data}).

\item[(4)]
Introduce the column vectors  $\mathbf{Y}=(Y_1,\ldots, Y_{n_2})^T$,
$\bar{\mathbf{x}}_k = ( \mathbf{x}_k^{\mathbf{i}}\dvtx \mathbf{i}\in
\mathbb{I}(r_\alpha,d))^T$, $k=1,\ldots,n_2$,
and form the data-matrix $\mathbf{X}=(\bar{\mathbf{x}}_1,
\ldots, \bar{\mathbf{x}}_{n_2})^T$ of
dimension $n_2\times (q(r_\alpha,d)+1)$.
Now using $\mathcal{D}_2^*$, fit a
polynomial regression model of order $r_\alpha$ by
\[
\tilde{\bolds{\vartheta}}=\arg\min_{\bolds{\vartheta}} \sum
_{k=1}^{n_2} \bigl(Y_k- f_{\bolds{\vartheta}}(
\mathbf{x}_k)\bigr)^2 =\arg\min_{\bolds{\vartheta}} \|
\mathbf{Y}-\mathbf{X} {\bolds{\vartheta}}\|^2,
\]
where the polynomial $f_{\bolds{\vartheta}}$ is introduced by
(\ref{multpolynomial}).
The unique least squares solution is given by
$\tilde{\bolds{\vartheta}} = (\mathbf{X}^T\mathbf{X})^{-1} \mathbf{X}^T \mathbf{Y}$,
since $\mathbf{X}$ is full-rank by Lemma~\ref{indeplemmult} below.
Intuitively, this is expected since the number of observations
$n_2\ge (2l+1)^d \ge (r_\alpha+1)^d \ge
|\mathbb{I}(r_\alpha,d)|= q(\alpha,d)+1$.

\item[(5)]
Finally, define the two-stage estimator $(\hat{\bolds{\mu}},\hat{M})$
of $(\bolds{\mu},M)$ by
\begin{equation}
\label{muestim1} \hat{\bolds{\mu}}=\arg\max_{\mathbf{x}\in C(\tilde{\bolds{\mu}},l\delta_n)}
f_{\tilde{\bolds{\vartheta}}}(\mathbf{x}), \qquad\hat{M}=f_{\tilde{\bolds{\vartheta}}}(\hat{\bolds{\mu}}).
\end{equation}
Note that $(\hat{\bolds{\mu}},\hat{M})$ depends on the first-stage
estimator and  the localization parameter $\delta_n$ introduced in
step 3.

\end{longlist}

Clearly the construction of the  two-stage procedure does not assume
the knowledge of the error variance $\sigma^2$ provided that the
preliminary estimator $\bolds{\tilde{\mu}}$ also does not use this
knowledge. Furthermore, the  two-stage approach simultaneously
estimates $\bolds{\mu}$ and $M$, since the same design points for both
estimators are used in the procedure. The two-stage procedure also
provides improved estimators for all the relevant derivatives of $f$
at $\bolds{\mu}$; see Remark~\ref{remlast} below.

The following theorem gives the rate of convergence of
the two-stage procedure for any smoothness level $\alpha>2$.

\begin{theorem}
\label{maintheorem}
Suppose that the localization parameter $\delta_n$
satisfies $\sqrt{n}\delta_n^2\to\infty$
and $\|\tilde{\bolds{\mu}}-\bolds{\mu}\| = o_p(\delta_n)$.
Then under conditions \textup{(A1)} and \textup{(A2)},
\begin{equation}
\label{new5} \|\hat{\bolds{\mu}}-\bolds{\mu}\| =O_p
\bigl(n^{-1/2}\delta_n^{-1} \bigr) +O_p
\bigl(\delta_n^{\alpha-1} \bigr)
\end{equation}
and
\begin{equation}
\label{new6} \hat{M} -M = O_p \bigl(n^{-1/2} \bigr)
+O_p \bigl(\delta_n^\alpha \bigr).
\end{equation}
\end{theorem}

Condition $\|\tilde{\bolds{\mu}}-\bolds{\mu}\|= o_p(\delta_n)$ has a clear
heuristic interpretation: at the second stage, one should not
localize more than what the accuracy of the estimation procedure
allows at the first stage.
Actually, it is sufficient to assume that $\P(\|\tilde{\bolds{\mu}}-\bolds{\mu}\|\le K\delta_n)\to 1$
for some $K$, but the dependence of $K$ on unknown quantities will complicate the analysis.

We first observe that there is always a rate improvement
from the first stage to the second if $\delta_n$ is chosen properly.
To see this, let $\varepsilon_n$ be the rate of convergence of~$\tilde{\bolds{\mu}}$. Since $\varepsilon_n$ cannot be better than the
optimal rate of convergence of all possible sequential procedures,
which is $n^{-(\alpha-1)/2\alpha}$, we have $\varepsilon_n\gtrsim
n^{-(\alpha-1)/(2\alpha)}$. Choose $\delta_n=\max(m_n \varepsilon_n,
n^{-1/(2\alpha)})$, where $m_n$ is a positive sequence going to
infinity sufficiently slowly. Then $\|\tilde{\bolds{\mu}}-\bolds{\mu}\|=
o_p(\delta_n)$ and  $\sqrt{n}\delta_n^2\to\infty$ (as $\alpha>2$) are
satisfied, and hence it remains to show
that $\delta_n^{\alpha-1}=O(\varepsilon_n)$ and $n^{-1/2}\delta_n^{-1}
=O(\varepsilon_n)$. If $\varepsilon_n\ll n^{-1/(2\alpha)}$,
then $\delta_n=n^{-1/(2\alpha)}$ and the second stage rate is
$n^{-1/2}\delta_n^{-1}=\delta_n^{\alpha-1}=n^{-(\alpha-1)/(2\alpha)}
=O(\varepsilon_n)$, and the order improves strictly unless
$\varepsilon_n\asymp n^{-(\alpha-1)/(2\alpha)}$. Clearly, in this case,
the choice of $\delta_n$ is optimal as it balances the ``order of
variability'' $n^{-1/2}\delta_n^{-1}$ and the ``order of bias''
$\delta_n^{\alpha-1}$. On the other hand, if $\varepsilon_n\gtrsim
n^{-1/(2\alpha)}$, $\delta_n=m_n \varepsilon_n$,
so $n^{-1/2}\delta_n^{-1}=o(\delta_n^{\alpha-1})$ and the second
stage rate is
$\delta_n^{\alpha-1}=m_n^{\alpha-1}\varepsilon_n^{\alpha-1}\ll
\varepsilon_n$, since $\alpha>2$ and $m_n$ grows sufficiently slowly.
Note that the ``optimal choice'' $\delta_n=n^{-1/(2\alpha)}$ is
prohibited in this case since we need $\varepsilon_n=o(\delta_n)$. For
estimating $M$, the rate of convergence of the
two-stage procedure clearly is $\max(n^{-1/2}, m_n^\alpha \varepsilon_n^\alpha)$,
which matches the optimal rate $n^{-1/2}$ if $\varepsilon_n\ll n^{-1/(2\alpha)}$.
Of course, if the choice of $\delta_n$ is too big, then the rates for
estimating $\bolds{\mu}$ or $M$ may deteriorate in the second stage.

Clearly, it is natural to use a preliminary estimator $\tilde{\bolds{\mu}}$ with the fastest possible
convergence rate $\varepsilon_n =n^{-(\alpha-1)/(2\alpha+d)}$ for any nonsequential procedure.
Then the two-stage estimator will lead to the best possible convergence rates
$n^{-(\alpha-1)/(2\alpha)}$ and $n^{-1/2}$ for estimating $\bolds{\mu}$ and $M$, respectively,
among all sequential procedures, provided that $\varepsilon_n = o(n^{-1/(2\alpha)})$.
This condition holds when $\frac{\alpha-1}{2\alpha+d} > \frac{1}{2\alpha}$,
or equivalently, $\alpha >1+\sqrt{1+d/2}$. Indeed, under this condition,
the two-stage procedure achieves the optimal rates $n^{-(\alpha-1)/(2\alpha)}$ and $n^{-1/2}$
for estimating $\bolds{\mu}$ and $M$, respectively, even when a rate-optimal estimator is not used,
as long as the convergence rate $\varepsilon_n$ of the preliminary estimator is faster than
$n^{-1/(2\alpha)}$.
If the condition $\varepsilon_n = o(n^{-1/(2\alpha)})$ fails, the two-stage procedure does not
give the optimal rate. In Section~\ref{sec4}, we discuss a multi-stage generalization that can
achieve optimal rate starting with almost any first-stage estimator.

The following corollary summarizes our conclusions.

\begin{corollary}
\label{corollary1}
Suppose that $\alpha>1+\sqrt{1+d/2}$ and conditions \textup{(A1)}, \textup{(A2)} hold.
If the convergence rate of the preliminary  estimator is faster than $n^{-1/(2\alpha)}$
and the localization parameter is $\delta_n=n^{-1/(2\alpha)}$, then
\[
\|\hat{\bolds{\mu}}-\bolds{\mu}\| = O_p \bigl(n^{-(\alpha-1)/2\alpha}
\bigr), \qquad\hat{M} -M = O_p \bigl(n^{-1/2} \bigr).
\]
\end{corollary}

Interestingly, dimension $d$, which affects the first-stage optimal convergence rates
$n^{-(\alpha-1)/(2\alpha+d)}$ and $n^{-\alpha/(2\alpha+d)}$ for estimating $\bolds{\mu}$ and $M$,
respectively, does not affect the corresponding two-stage and fully sequential optimal
convergence rates $n^{-(\alpha-1)/(2\alpha)}$ and $n^{-1/2}$. Thus the curse of dimensionality
is nearly avoided by the two-stage procedure,
provided that the regression function is sufficiently smooth to ensure
$\alpha >1+\sqrt{1+d/2}$. The lower bound in this inequality
increases with the dimension $d$.
Notice that if $\alpha >3$, the corollary yields the optimal rates
for estimating $\bolds{\mu}$ and $M$ for all the dimensions for which $3\ge
1+\sqrt{1+d/2}$, that is, up to dimension $d=6$, including  the most
important dimensions $d=1,2,3$.

\begin{remark}
\label{remark1}
We can formulate a uniform version of Theorem~\ref{maintheorem}.
By inspecting the proofs, we see that all the bounds for the two-stage
procedure can be made uniform over the H\"older class\vadjust{\goodbreak}
$\mathcal{H}_d(\alpha,L, D)$ if we additionally require
relation (\ref{newA3}) for some $\kappa,\lambda_0>0$,
the uniform boundedness of all the partial derivatives involved in
the definition of $\mathcal{H}_d$ and the uniformity
of the first stage estimator.

To be more specific, for some positive $\alpha$, $L$, $L_1$,
$\kappa$, $\lambda_0$, $\kappa_1$, $\delta$ and $\varepsilon$, such that
\mbox{$\alpha>2$} and $\kappa_1\le\kappa$, and a compact convex
$D\subseteq \mathbb{R}^d$, introduce the following conditions:

\begin{longlist}[(\~A1)]
\item[(\~A1)]
$f\in \mathcal{H}_d(\alpha,L, D^\varepsilon)$ and
$\sup_{\mathbf{x} \in D^\varepsilon}  |\mathbf{D}^{\mathbf{i}} f (\mathbf{x}) | \le L_1$
for all $\mathbf{i} \in \mathbb{I}(r_\alpha,d)$.

\item[(\~A2)]
There is a unique point $\bolds{\mu}\in\mathring{D}$
that maximizes  the function $f$ on $D$,
$\sup_{\mathbf{x}\in D} f(\mathbf{x})
=\max_{\mathbf{x}\in \mathring{D}} f(\mathbf{x}) = f(\bolds{\mu})$.
Moreover, $f(\bolds{\mu})\ge f(\mathbf{x}) + \delta$ for all $\mathbf{x}
\notin B(\bolds{\mu}, \kappa_1)$ and
$\sup_{\mathbf{x}\in B(\bolds{\mu}, \kappa)}
\lambda_{\max}(Hf(\mathbf{x}))\le -\lambda_0$.
\end{longlist}
Let $\tilde{\mathcal{H}}_d$ be the class
of functions which satisfy (\~A1) and (\~A2).
Then Theorem~\ref{maintheorem} holds uniformly in $f \in \tilde{\mathcal{H}}$,
provided $\|\tilde{\bolds{\mu}}-\bolds{\mu}\| = o_p(\delta_n)$ holds
uniformly over $\tilde{\mathcal{H}}$.
Condition (\~A1) is a strengthened version of (A1), namely (A1) is complemented
by the requirement of uniform boundedness  of all the relevant partial derivatives.
Condition (\~A2) is in turn a stronger version of (A2): relation (\ref{newA3})
is included in (\~A2) with common $\kappa$ and $\lambda$ for the whole class, and
the existence of a unique location~$\bolds{\mu}$ of maximum
is strenthened by the requirement of the uniform separation of
the maximum function value $f(\bolds{\mu})$
from the function values outside $B(\bolds{\mu}, \kappa_1)$.
Inside this vicinity, as $\kappa_1\le\kappa$, the separation of the maximum
can be characterized by the Taylor expansion and (\ref{newA3});
see the arguments in (\ref{lem1rel1}) below. This uniform
separation condition is essential  to make the first stage rate for
$\tilde{\bolds{\mu}}$ uniform over the functional class. Note that the separation condition
for any particular function holds by the compactness of $D$ and the
uniqueness of the location of the maximum.

On the other hand, the two-stage procedure can achieve improved rates only
under a local H\"older condition satisfied in a neighborhood of $\bolds{\mu}$
provided that a first stage estimator with sufficiently good rate is available
as a preliminary estimator. This is due to the fact that the second stage design
points are chosen close to the preliminary estimate, and hence close to
the true maximum location $\bolds{\mu}$.
\end{remark}

\begin{remark}
Almost sure convergence of $\hat{\bolds{\mu}}$ and $\hat{M}$ can be obtained
assuming that the preliminary estimator convergence rate is given in the almost sure sense.
This will follow from the estimates given in Lemmas~\ref{apprlemmult} and~\ref{Zlemmult}. Under additional moment conditions on the error distribution,
almost sure convergence rate of a kernel-type estimator can be found.
\end{remark}

\section{Multi-stage procedures and resolving the curse of dimensionality}\label{sec4}
Theorem~\ref{maintheorem} shows that for estimating the maxima
$\bolds{\mu}$ and maximum value $M$ of a H\"older $\alpha$-smooth
function $f$, $\alpha>2$, starting
with an estimator $\tilde{\bolds{\mu}}$ having convergence rate
$\varepsilon_n$ and localization parameter $\delta_n$, a two-stage
estimator has an improved rate of convergence, unless
$\varepsilon_n$ is already equal to the optimal rate $n^{-(\alpha-1)/(2\alpha)}$. More precisely,
assuming that $\varepsilon_n=O(n^{-\gamma})$ for some $\gamma>0$ and choosing
$\delta_n=\max(m_n\varepsilon_n, n^{-1/(2\alpha)})$, where $m_n\to\infty$ is a slowly
varying sequence, the convergence rate for estimating $\bolds{\mu}$ improves to
the optimal rate $n^{-(\alpha-1)/(2\alpha)}$ if $\varepsilon_n=o(n^{-1/(2\alpha)})$ and to
$\varepsilon_n^{\alpha-1}$ up to a slowly varying factor, if $\varepsilon_n\gtrsim n^{-1/(2\alpha)}$.
Although the latter rate is not optimal, further improvement in
rate can be achieved by applying the two-stage technique again, using
the estimator obtained in the second stage as
the new preliminary estimator, and repeating the procedure until
the optimal rate is obtained. After $k$ iterations of the two-stage
procedure, the convergence rate thus becomes
$\varepsilon_n^{(\alpha-1)^k}$ up to a slowly varying factor, provided
that $\varepsilon_n^{(\alpha-1)^k}\gg n^{-(\alpha-1)/(2\alpha)}$.
Let $k_1$ be the largest integer such that the last relation holds for
$k=k_1$. Then iterating the two-stage procedure $k_1+1$ times, the
resulting estimator will have the optimal convergence rate
$n^{-(\alpha-1)/(2\alpha)}$. Thus the final multi-stage procedure has
convergence rate completely free of the dimension $d$ and applies to
any smoothness level $\alpha>2$. In order to apply the procedure in
$k_1+1$ stages, one will need to split the observation
budget in $k_1+1$ parts following the description given in step 3 of
the procedure.

If we are interested only in estimating the maximum $M$, we may be
able to stop earlier when applying the multi-stage procedure. In this
case, the target optimal rate is $n^{-1/2}$. The two-stage estimator
has convergence rate given by $\max\{\varepsilon_n^\alpha, n^{-1/2}\}$.
Hence the optimal rate will be obtained at stage $k_2+1$, where
$k_2$ is the largest integer integer such that
$\varepsilon_n^{\alpha^k}\gg n^{-1/2}$.

\begin{remark}
The smoothness level $\alpha$ needs to be strictly greater than
$2$ to control the error in the second-order Taylor approximation of the underlying
multivariate regression function. As $\alpha$ gets closer to $2$, the required
number of stages in the multi-stage procedure increases without bound.
\end{remark}

Consider now the adaptive version of our original estimation problem,
where the problem is to estimate $\bolds{\mu}$ at the optimal rate
$n^{-(\alpha-1)/(2\alpha)}$ and $M$ at rate $n^{-1/2}$ without knowing
the smoothness level $\alpha$. Since the choice of the localization
parameter $\delta_n$ depends on the knowledge of $\alpha$ it is not
possible to apply the two-stage procedure, and hence a multi-stage
adaptive estimator for $\bolds{\mu}$ is not possible. However, for
estimating $M$, it is possible to construct a multi-stage procedure
with convergence rate $n^{-1/2}$ without knowing $\alpha$, as long as
$\alpha>2$. Start with an estimator for $\bolds{\mu}$ which converges at
rate $n^{-\beta}$ for all H\"older 2-smooth functions. For instance,
the rate $n^{-1/(4+d)}$ is possible in dimension $d$ by the results
of M\"uller (\citeyear{Muller1989}) and Facer and M\"uller (\citeyear{FacerMuller2003}). Then by applying
Theorem~\ref{maintheorem} with $\delta_n=m_n n^{-\beta}$, where
$m_n\to\infty$ is a slowly varying sequence, the rate of convergence
for estimating $M$ improves to $\max\{m_n^2 n^{-2\beta}, n^{-1/2}\}$
in stage two. Repeating the two-stage procedure $k$ times, thus
the rate will improve to $n^{-1/2}$, whenever $2^k\beta> \frac{1}{2}$,
or $k> (\log(1/\beta)/\log 2)-1$. In particular, starting with the
one-stage optimal estimator having convergence rate $n^{-1/(4+d)}$,
the required number of stages to achieve $n^{-1/2}$ rate at all
H\"older 2-smooth functions is the smallest integer greater than
$(\log(4+d)/\log 2)-1$, since  repeating the two-stage procedure
beyond $k$ given above does not hurt the rate.

\begin{figure}

\includegraphics{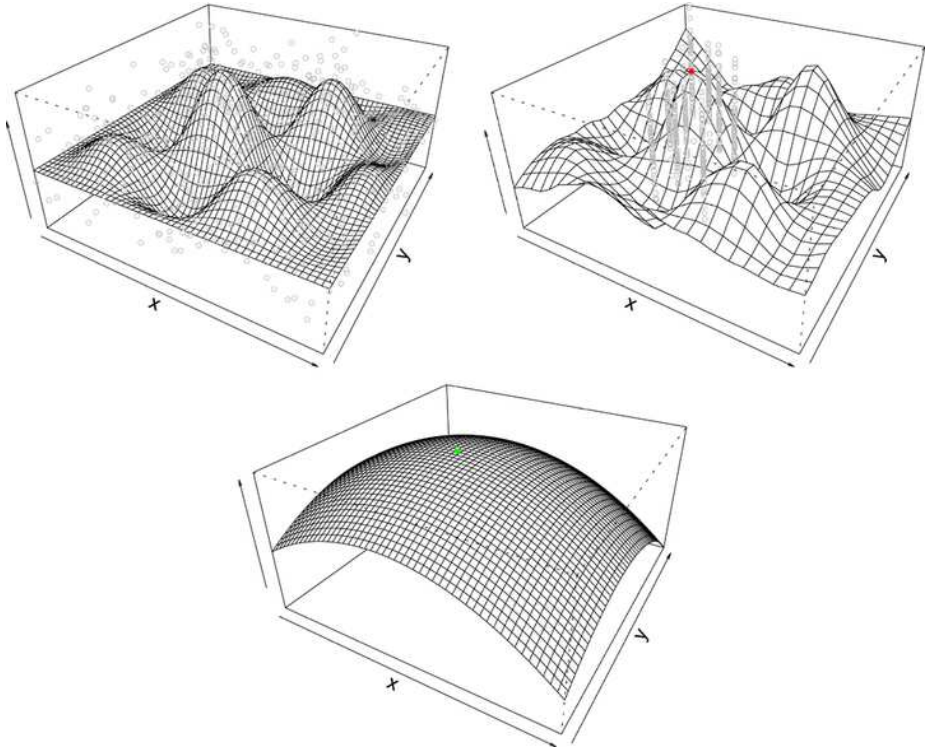}

\caption{Top left: true regression function (surface) and stage one observations
(gray points). Top right: surface fitted through stage one observations,
the initial estimator (red point) and the stage two observations (gray points).
Bottom: quadratic surface fitted through the stage two observations and
final estimator (green point).}
\label{fig3d1}
\end{figure}

\section{Simulations}\label{sec5}

In this section, we compare the performance of the two-stage procedure with
an equivalent single-stage procedure. We consider the bivariate case $d=2$ and
take a regression function $f\dvtx [0,1]^2 \to \mathbb{R}$ defined by
\[
f(x,y) = 5x(x-1)y(y-1)\sin(11x)\sin(11y)
\]
(the smooth surface in the top left panel of Figure~\ref{fig3d1}).
In the first stage the function is observed with Gaussian noise with standard
deviation $\sigma=0.1$ on a regular $25$ by $25$ grid (gray points in the same panel).
Using standard local linear regression, a surface is fit through these points
(the surface in the top right panel of Figure~\ref{fig3d1})
and the point where this fitted function is maximal serves as the stage one
estimator $\tilde{\bolds{\mu}}$ (the red point in the same panel).
Next we take $\delta = 0.1$ and generate $70$ new observations at each
of the nine points  $(\tilde\mu_1+j_1\delta, \tilde\mu_2 + j_2\delta)$,
$j_1, j_2 =0, \pm 1$ (gray points in the top right panel of Figure~\ref{fig3d1}).
Finally a quadratic surface is fitted through these new data points (the surface in
the lower panel of Figure~\ref{fig3d1}) and  the location of the maximum is
the final second stage estimator $\hat{\bolds{\mu}}$ (the green point in the figure).
The implementation of the procedure is rather straightforward. In the statistical
language R, we used the standard function \texttt{loess} in the first stage to
fit the surface using the first stage observations, and we used the function
\texttt{lm} to fit the quadratic surface using the second stage observations.

\begin{figure}[b]

\includegraphics{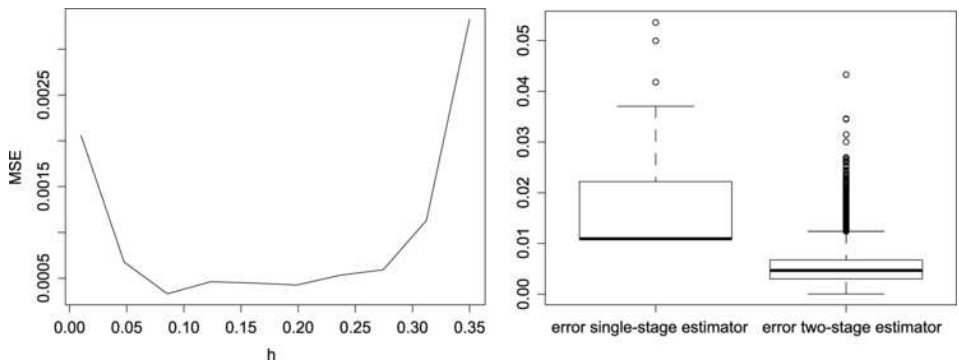}

\caption{Left: MSE of the single-stage estimator against the bandwidth used.
Right: boxplot of the errors of the single-stage estimator
(with optimal bandwidth choice $h = 0.085$) and the errors of our two-stage
procedure (with the same bandwidth choice and $\delta =0.1$).}
\label{fig3d2}
\end{figure}

Note that in total we have used $25\times 25 + 9\times 70 =  1255$ observations.
It is illustrative to compare our procedure to a single stage estimator that uses about the
same amount of regularly spaced observations. The closest is a regular
$36$ by $36$ grid, which contains $1296$ points. We make noisy observations
of the function~$f$ at these grid points, again corrupted by centered
Gaussian noise with standard deviation $0.1$.
We consider the estimator for the location of the maximum of $f$ that is obtained
by fitting a locally linear surface through these data points and computing the
location where this is maximal.
Obviously, the quality of this estimator depends on the bandwidth that is used
(or span parameter, as it is called in the R function \texttt{loess}).
To obtain a fair comparison with our two-stage estimator, we should make
an optimal choice. We achieve this by repeating the experiment a large number
of times with different bandwidths and computing numerical mean squared errors (MSEs).
The result is shown in the left panel of Figure~\ref{fig3d2}.
The numerical MSE is minimal for the bandwidth choice $h = 0.085$.

To compare the mean-squared error of the single-stage estimator based on
this regular grid, we replicated the experiment $10\mbox{,}000$ times and computed
the Monte-Carlo average of the squared difference between the estimate and
the true maximum. The results are shown in the left boxplot in the right panel
of Figure~\ref{fig3d2}. Similarly we carried out the two-stage procedure
$10\mbox{,}000$ times (with bandwidth $h=0.085$ in stage one
and $\delta=0.1$) and computed the errors as well.
These are shown in the right boxplot in the right panel of Figure~\ref{fig3d2}.
It is clear that the two-stage estimator performs better in this situation,
in terms of the mean-squared error. This is in spite of the fact that the two-stage estimator
has used less observations, namely $1255$ in total compared to $1296$
used by the single-stage estimator.

\begin{figure}

\includegraphics{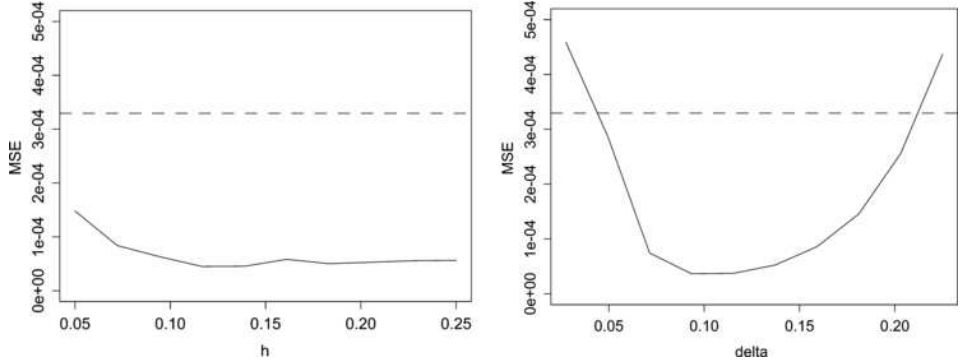}

\caption{MSE of the two-stage estimator against the bandwidth $h$ used
in the first stage (left), and against the localization
parameter $\delta$ (right). The dashed line is the MSE of the single-stage estimator
with optimal bandwidth choice.}
\label{fig3d3}
\end{figure}

In practice the quality of our procedure clearly depends on the quality of the estimator
that is used in the first stage and also on the choice of the localization parameter $\delta$.
In this simulation example, where we use local linear regression in the first stage,
the quality of the  estimator therefore depends on the bandwidth used in  stage one.
To investigate the dependence of the performance on this parameter
we carried out the simulation study described above for a range of bandwidths.
The results are shown in the left panel of Figure~\ref{fig3d3}.
The solid line gives the MSE of our estimator as a function of the bandwidth
used in stage one. The dashed line is the MSE of the optimal single stage estimator
described above. The plot shows that in fact for a range
of bandwidths the two stage procedure performs better than the single stage procedure.
Similarly, the right panel of Figure~\ref{fig3d3} describes the performance
of the two-stage procedure as a function of the localization parameter~$\delta$.
Again there is a range of possible $\delta$'s for which we
obtain an improved performance, but choosing $\delta$ too small
or too large deteriorates the quality. In practice one might, for instance,
use cross-validation type methods to set the tuning parameters $h$ and
$\delta$. Further research is needed to find theoretically sound methods.

\section{Proofs}\label{sec6}

Throughout this section, $\alpha$ and $d$ are are kept fixed.
To simplify notation, we abbreviate $\mathbb{I}(r_\alpha,d)$
by $\mathbb{I}$, $\mathbb{I}_k(d)$ by $\mathbb{I}_k$ and $q(\alpha,d)$ by $q$.

First we introduce several quantities we are going to use in the sequel.
Define $\mathbf{z}_k = \mathbf{x}_k - \tilde{\bolds{\mu}}$,
$k=1,\ldots, n_2$,
and reformulate definition (\ref{muestim1}) by
representing the involved quantities
in terms of the newly defined shifted design points
$\mathbf{z}_k$,  $k=1,\ldots, n_2$.
Let $\mathbf{d}_j= \tilde{\mathbf{d}}_j-\tilde{\bolds{\mu}}$,
$j=1,\ldots, (2l+1)^d$.
Then for all $k=1,\ldots, n_2$,
\begin{equation}
\label{designmult} \mathbf{z}_k \in \{0, \pm
\delta_n, \ldots, \pm l\delta_n\}^d = \{
\mathbf{d}_1,\ldots, \mathbf{d}_{(2l+1)^d}\}\subset C(l
\delta_n),
\end{equation}
so that each of the distinct $(2l+1)^d$ points are repeated $n_3=n_2/(2l+1)^d$
times in the new design set $\{\mathbf{z}_k,   k=1,\ldots, n_2\}$.
Using definition (\ref{multpolynomial}), define an estimator~$\hat{\bolds{\theta}}$ by equating the two polynomials
%
\begin{equation}
\label{defhattheta} f_{\hat{\bolds{\theta}}}(\mathbf{x}-
\tilde{\bolds{\mu}}) =f_{\tilde{\bolds{\vartheta}}}(\mathbf{x}) \qquad\mbox{equivalently } \hat{
\bolds{\theta}} = \bigl(\mathbf{Z}^T\mathbf{Z}\bigr)^{-1}
\mathbf{Z}^T \mathbf{Y},
\end{equation}
where
\begin{equation}
\label{defZmult} \mathbf{Z}=(\bar{
\mathbf{z}}_1, \ldots, \bar{\mathbf{z}}_{n_2})^T,\qquad
\bar{\mathbf{z}}_k = \bigl(\mathbf{z}_k^{\mathbf{i}}\dvtx \mathbf{i}\in \mathbb{I}\bigr)^T
\end{equation}
and $\mathbf{z}_k = \mathbf{x}_k - \tilde{\bolds{\mu}}$,
$k=1,\ldots,n_2$.
The matrix $\mathbf{Z}^T\mathbf{Z}$ is invertible by Lemma
\ref{indeplemmult} below.
We thus obtain an equivalent description of
the estimator $(\hat{\bolds{\mu}}, \hat{M})$
given by (\ref{muestim1})
in terms of the polynomial $f_{\hat{\bolds{\theta}}}(\mathbf{z})$
defined by (\ref{multpolynomial}),
with $\hat{\bolds{\theta}}$ defined by (\ref{defhattheta}),
\begin{equation}
\label{muestim2} \hat{\bolds{\mu}} =\tilde{\bolds{\mu}} +
\mathring{\bolds{\mu}},\qquad \hat{M} = f_{\hat{\bolds{\theta}}}(\mathring{\bolds{\mu}})
\qquad\mbox{where } \mathring{\bolds{\mu}}=\arg\max_{\mathbf{z}\in C(l\delta_n)} f_{\hat{\bolds{\theta}}}(
\mathbf{z}),
\end{equation}
$C(l\delta_n) =[- l\delta_n,l\delta_n]^d\subset \mathbb{R}^d$.
In doing this shifting trick, we make the computations easier
because the matrix $\mathbf{Z}^T\mathbf{Z}$ will have a lot of zero entries as
the design points $\mathbf{z}_k$'s are symmetrically
centered around zero in each dimension rather than being centered around $\tilde{\bolds{\mu}}$.

Next, let the vector $\bolds{\theta} = (\theta_{\mathbf{i}}\dvtx \mathbf{i}\in \mathbb{I})^T$ be defined
by the equality of the two polynomials $f_{\bolds{\theta}}(\mathbf{x}-\tilde{\bolds{\mu}})
= P_{f,\bolds{\mu}}(\mathbf{x})$, where $P_{f,\bolds{\mu}}$ is the Taylor
expansion of $f$ of order $r_\alpha$ around $\bolds{\mu}$ defined by
(\ref{new3}),
\begin{equation}
\label{multdeftheta} \sum
_{\mathbf{i} \in \mathbb{I}} \theta_{\mathbf{i}}(\mathbf{x}-\tilde{\bolds{
\mu}})^{\mathbf{i}} =f(\bolds{\mu})+\sum_{\mathbf{i} \in \mathbb{I}, |\mathbf{i}|\ge 2}
\frac{\mathbf{D}^{\mathbf{i}} f(\bolds{\mu})}{\mathbf{i}!}
(\mathbf{x}-\bolds{\mu})^{\mathbf{i}},
\end{equation}
here we have used the condition $\nabla f(\bolds{\mu})= \mathbf{0}$, due to (A1) and (A2).
Thus, $\bolds{\theta}$ is a random vector depending on $f$, $\bolds{\mu}$ and $\tilde{\bolds{\mu}}$.
From (\ref{multdeftheta}) it follows that
\begin{equation}
\label{new9} \mathbf{i}! \theta_{\mathbf{i}} = \mathbf{D}^{\mathbf{i}}
P_{f,\bolds{\mu}}(\tilde{\bolds{\mu}}),\qquad  \mathbf{D}^{\mathbf{i}} f(\bolds{\mu})
=\mathbf{D}^{\mathbf{i}} f_{\bolds{\theta}}(\bolds{\mu}-\tilde{\bolds{\mu}}),\qquad
\mathbf{i} \in \mathbb{I}.
\end{equation}

The next lemma ensures that the estimator  (\ref{defhattheta}) is
well defined;  that is, $\mathbf{Z}^T\mathbf{Z}$ is invertible.

\begin{lemma}
\label{indeplemmult}
The columns of matrix $\mathbf{Z}$ (and $\mathbf{X}$)
defined by (\ref{defZmult})
are linearly independent.
\end{lemma}

\begin{pf}
Consider the matrix $\mathbf{Z}$;
the same proof applies to $\mathbf{X}$.

For multi-indices $\mathbf{i}\in \mathbb{N}^p$ and
$\mathbf{j}\in \mathbb{N}^s$, define
the concatenation operation
$(\mathbf{i}, \mathbf{j})=(i_1,\ldots, i_p, j_1,\ldots, j_s)
\in \mathbb{N}^{p+s}$. In particular,
for $k\in \mathbb{N}$, $\mathbf{i} \in \mathbb{N}^p$,
$(k,\mathbf{i}) = (k, i_1,\ldots, i_p)$.
Introduce the following notation: for $l=0,1,\ldots, d-1$
and $\mathbf{z}=(z_1,\ldots,z_d)\in \mathbb{R}^d$, define
$\mathbf{i}_{-l} = (i_{l+1}, \ldots, i_d) \in \mathbb{N}^{d-l}$,
$\mathbf{z}_{-l} = (z_{l+1}, \ldots, z_d) \in \mathbb{R}^{d-l}$ and the set
$\mathbb{I}_{-l}(i_1,\ldots, i_l)=
\{\mathbf{i} \in \mathbb{N}^{d-l}\dvtx (i_1,\ldots ,i_l,\mathbf{i}) \in \mathbb{I}\}$.

Let  $\mathbf{C}_{\mathbf{i}}$, $\mathbf{i}\in\mathbb{I}$,
be the columns of the matrix $\mathbf{Z}$. We need to show that
$\sum_{\mathbf{i} \in \mathbb{I}}
\lambda_{\mathbf{i}} \mathbf{C}_{\mathbf{i}} =\mathbf{0}$ implies that
$\lambda_{\mathbf{i}}=0$ for all $\mathbf{i}\in\mathbb{I}$.
The equality $\sum_{\mathbf{i} \in \mathbb{I}}
\lambda_{\mathbf{i}} \mathbf{C}_{\mathbf{i}} =\mathbf{0}$ is
equivalent to
\[
0 = \sum_{\mathbf{i} \in \mathbb{I}} \lambda_{\mathbf{i}}
\mathbf{z}_k^{\mathbf{i}}, \qquad k=1,2,\ldots, n_2.
\]

Among $\{\mathbf{z}_1,\ldots, \mathbf{z}_{n_2}\}$,
only $(2l+1)^d$ are distinct---$\{\mathbf{d}_1,\ldots, \mathbf{d}_{(2l+1)^d}\}$ given by~(\ref{designmult}).
Thus, for all
$\mathbf{z}\in\{\mathbf{d}_1,\ldots, \mathbf{d}_{(2l+1)^d}\}$,
\[
0 = \sum_{\mathbf{i} \in \mathbb{I}} \lambda_{\mathbf{i}}
\mathbf{z}^{\mathbf{i}} = \sum_{i_1=0}^{r_\alpha}
z_1^{i_1} \sum_{\mathbf{i}_{-1}\in \mathbb{I}_{-1}(i_1)}
\lambda_{i_1 \mathbf{i}_{-1}} \mathbf{z}_{-1}^{\mathbf{i}_{-1}}.
\]
For a fixed $\mathbf{z}_{-1}=(z_2,\ldots,z_d)$,
the right-hand side of the last relation
is a polynomial of order $r_\alpha$ in variable $z_1$.
But we have $2l+1>r_\alpha$ different design values $\{j\delta_n\dvtx j=0,\pm1, \pm 2, \ldots, \pm l\}$
of the variable $z_1$ for which this polynomial
must take the zero value. This forces all the coefficients of this polynomial to be
zero. Thus we have that
\[
0=\sum_{\mathbf{i}_{-1}\in \mathbb{I}_{-1}(i_1)} \lambda_{i_1 \mathbf{i}_{-1}}
\mathbf{z}_{-1}^{\mathbf{i}_{-1}},\qquad  i_1=0,1,\ldots,
r_\alpha
\]
for all possible design values of $\mathbf{z}_{-1}=(z_2,\ldots, z_d)$.
Iterating the above reasoning up to the variable $z_d$ leads to,
for all $i_1,\ldots,i_{d-1}=0,1,\ldots,r_\alpha$,
$z_d\in\{0,\pm \delta_n, \pm 2 \delta_n, \ldots, \pm l\delta_n\}$,
\[
0 = \sum_{i_d\in \mathbb{I}_{-(d-1)}(i_1,i_2,\ldots, i_{d-1})} \lambda_{i_1i_2, \ldots, i_{d-1} i_d}
z_d^{i_d},
\]
from which we derive that $\lambda_{\mathbf{i}}=0$ for all $\mathbf{i}\in\mathbb{I}$.
\end{pf}

\begin{remark}
In the case $d=1$,  $\mathbf{X}$ and $\mathbf{Z}$ are Vandermonde matrices.
\end{remark}

The next lemma shows that the second stage data $\mathcal{D}_2^*$
can be regarded as coming approximately from a certain polynomial regression model.

\begin{lemma}
\label{apprlemmult}
Assume \textup{(A1)} and let the data $\{(\mathbf{x}_k,Y_k),   k=1,\ldots, n_2\}$
and $\bolds{\xi}= (\xi_1,\ldots, \xi_{n_2})^T$
be given by the second stage observation scheme, $\mathbf{Y}=(Y_1,\ldots, Y_{n_2})^T$,
$\bolds{\eta} =(\eta_1,\ldots, \eta_{n_2})^T$ with
$\eta_k=f(\mathbf{x}_k) - P_{f,\bolds{\mu}}(\mathbf{x}_k)$. Then\vadjust{\goodbreak}
$\mathbf{Y}= \mathbf{Z}{\bolds{\theta}} +\bolds{\eta}+\bolds{\xi}$,
where $\mathbf{Z}$ and  ${\bolds{\theta}}$ are defined
by (\ref{defZmult}) and (\ref{multdeftheta}), respectively, and
$\bolds{\eta}$ is independent of $\bolds{\xi}$.
Moreover,  for some constants $C_1, C_2$ and uniformly in $k \in \{1,2,\ldots,n_2\}$,
\begin{equation}
\label{multeta} |\eta_k| \le C_1
\delta_n^\alpha + C_2 \|\tilde{\bolds{\mu}}-
\bolds{\mu}\|^\alpha.
\end{equation}
\end{lemma}

\begin{pf}
Since  $\eta_k=f(\mathbf{x}_k) - P_{f,\bolds{\mu}}(\mathbf{x}_k)$,
by (\ref{defZmult}) and (\ref{multdeftheta}),
\begin{eqnarray*}
Y_k &=& f(\mathbf{x}_k) + \xi_k =
P_{f,\bolds{\mu}}(\mathbf{x}_k)+\eta_k+\xi_k
\\
&=& f_{{\bolds{\theta}}}(\mathbf{x}_k-\tilde{\bolds{\mu}}) +
\eta_k +\xi_k = \bar{\mathbf{z}}_k^T
{\bolds{\theta}} +\eta_k +\xi_k.
\end{eqnarray*}
Clearly, $\bolds{\eta}$ is independent of $\bolds{\xi}$ by definition.
It remains to show (\ref{multeta}).
Apply the $c_r$-inequality, $|a+b|^r \le \max(1,2^{r-1}) (|a|^r + |b|^r)$, $r>0$,
(\ref{multHolder})  and the fact that
$\|\mathbf{x}_k-\tilde{\bolds{\mu}}\|\le \sqrt{d}l\delta_n$, $k=1,\ldots, n_2$,
to obtain (\ref{multeta})
\begin{eqnarray*}
|\eta_k|&=& \bigl|f(\mathbf{x}_k) - P_{f,\bolds{\mu}}(
\mathbf{x}_k)\bigr| \le L \|\mathbf{x}_k-\bolds{\mu}
\|^\alpha
\\[-1pt]
&\le& L c_\alpha \bigl(\|\mathbf{x}_k-\tilde{\bolds{\mu}}
\|^\alpha + \|\tilde{\bolds{\mu}} -\bolds{\mu}\|^\alpha \bigr) \le
C_1 \delta_n^\alpha + C_2 \|\tilde{
\bolds{\mu}}-\bolds{\mu}\|^\alpha.\vspace*{-1pt}
\end{eqnarray*}
\upqed\end{pf}

Lemma~\ref{indeplemmult} ensures that  the matrix $\mathbf{Z}^T\mathbf{Z}$
is nonsingular. The following lemma describes the asymptotic behavior of the
elements of its inverse. For notational convenience, below we enumerate the
rows and columns of matrices starting from~$0$.

Enumerate $\mathbb{I}$ by arranging their elements in the
order described in Section~\ref{sec2}, which we denote by $\mathbf{i}_0,\ldots,\mathbf{i}_q$, respectively.

\begin{lemma}
\label{Zlemmult}
The $(i,j)$th element $h_{ij}$ of $(\mathbf{Z}^T\mathbf{Z})^{-1}$ satisfies
\[
h_{ij} =O \bigl(n^{-1} \delta_n^{-(|\mathbf{i}_i|+|\mathbf{i}_j|)}
\bigr),\qquad  i,j =0,1,\ldots,q.
\]
\end{lemma}

\begin{pf}
Since $\mathbf{z}_k^{\mathbf{0}}=1$
for all $k=1,\ldots, n_2$, we have
\[
\mathbf{Z}^T\mathbf{Z}=\pmatrix{ n_2 & \displaystyle\sum
_{k=1}^{n_2} \mathbf{z}_k^{\mathbf{i}_1} &
\cdots & \displaystyle\sum_{k=1}^{n_2}
\mathbf{z}_k^{\mathbf{i}_q} \vspace*{2pt}
\cr
\displaystyle\sum
_{k=1}^{n_2} \mathbf{z}_k^{\mathbf{i}_1} &
\displaystyle\sum_{k=1}^{n_2} \mathbf{z}_k^{\mathbf{i}_1}
\mathbf{z}_k^{\mathbf{i}_1} &\cdots & \displaystyle\sum
_{i=1}^{n_2} \mathbf{z}_k^{\mathbf{i}_1}
\mathbf{z}_k^{\mathbf{i}_q}\vspace*{2pt}
\cr
&\cdots& \vspace*{2pt}
\cr
\displaystyle\sum_{k=1}^{n_2} \mathbf{z}_k^{\mathbf{i}_q}
& \displaystyle\sum_{i=1}^{n_2} \mathbf{z}_k^{\mathbf{i}_1}
\mathbf{z}_k^{\mathbf{i}_q}&\cdots & \displaystyle\sum
_{i=1}^{n_2} \mathbf{z}_k^{\mathbf{i}_q}
\mathbf{z}_k^{\mathbf{i}_q} }.
\]
Then, for some constants $a_{ij}$, $i,j=0,1,\ldots, q$, we rewrite
the symmetric matrix $\mathbf{Z}^T\mathbf{Z}$ as follows:\vspace*{-1pt}
\[
\mathbf{Z}^T\mathbf{Z} = n_3 \pmatrix{
a_{00} & a_{01} \delta_n^{|\mathbf{i}_1|} &
\cdots & a_{0q} \delta_n^{|\mathbf{i}_q|} \vspace*{2pt}
\cr
a_{10} \delta_n^{|\mathbf{i}_1|} & a_{11}
\delta_n^{|\mathbf{i}_1|+|\mathbf{i}_1|} &\cdots & a_{1q}
\delta_n^{|\mathbf{i}_1|+|{\mathbf{i}_q}| }
\vspace*{2pt}\cr
&\cdots& \vspace*{2pt}
\cr
a_{q0} \delta_n^{|\mathbf{i}_q|}
& a_{q1} \delta_n^{|\mathbf{i}_q|+|\mathbf{i}_1|} &\cdots &
a_{qq} \delta_n^{|\mathbf{i}_q|+|\mathbf{i}_q|}}.
\]
Some entries are easy to compute. For example,
$a_{00} = (2l+1)^d$ since $n_2=(2l+1)^d n_3$. Moreover,
there are many zeros due to the symmetry of the design.
In particular, $a_{ij}=0$ for all $i,j \in \{0,1,\ldots, q\}$ such that
$|\mathbf{i}_i|+|{\mathbf{i}_j}| $ is an odd number.
However we are not concerned about the exact values $a_{ij}$ but only about
nonsingularity of the matrix $\mathbf{A}$ with $(i,j)$th entry equal to $a_{ij}$, $i,j=0,\ldots,q$.

Let $\bolds{\Delta}$
be the diagonal matrix with elements
$\delta_n^{|\mathbf{i}_j|}$, $j=0,1,\ldots, q$, in that order. Now notice that
$\mathbf{Z}^T\mathbf{Z} = n_3\bolds{\Delta} \mathbf{A} \bolds{\Delta}$.
Since $\mathbf{Z}^T\mathbf{Z}$ is nonsingular
by Lemma~\ref{indeplemmult},
it follows that $\mathbf{A}$ is also invertible.
Therefore $ (\mathbf{Z}^T\mathbf{Z} )^{-1}
= n_3^{-1}\bolds{\Delta}^{-1} \mathbf{A}^{-1} \bolds{\Delta}^{-1}$.
Denote by ${a}^{ij}$ the $(i,j)$th entry of the constant matrix
$\mathbf{A}^{-1}$ and recall that $n_3 \ge c n$.
Then for $i,j =0,1,\ldots,q$,
\[
h_{ij} = n_3^{-1} {a}^{ij}
\delta_n^{-(|\mathbf{i}_i|+|\mathbf{i}_j|)}= O \bigl(n^{-1}
\delta_n^{-(|\mathbf{i}_i|+|\mathbf{i}_j|)} \bigr).
\]
\upqed\end{pf}

\begin{remark}
For $d=1$ and even $r_\alpha$, we have $2l+1=r_\alpha+1$.
Put  $b_0=r_\alpha+1$, $b_m= 0$ for all odd $m \in \{1,\ldots, 2r_\alpha\}$,
and for each even $m \in \{1,\ldots, 2r_\alpha\}$
\[
b_m = 2\bigl(1+2^m+3^m+\cdots +
l^m\bigr)=2\bigl\{1+2^m+\cdots + (r_\alpha/2)^m
\bigr\}.
\]
Then the entries of
$\mathbf{A}$ can be computed as follows:
Since $n_2=(2l+1)n_3$, $\sum_{k=1}^{n_2} z_k^m = 0$ for each odd
$m \in \{1,\ldots, 2r_\alpha\}$ and
\[
\sum_{k=1}^{n_2} z_k^m
= 2n_3\bigl\{l^m \delta_n^m
+(l-1)^m\delta_n^m + \cdots +
\delta_n^m\bigr\} = n_3\delta_n^m
b_m
\]
for each even $m \in \{1,\ldots, 2r_\alpha\}$, we obtain that $a_{ij} = b_{i+j}$.

The case of odd $r_\alpha$ can be treated similarly leading to slightly
different constants.
\end{remark}

\begin{lemma}
\label{thetalemmult}
Assume \textup{(A1)}, and let  $\hat{{\bolds{\theta}}}$ and  ${\bolds{\theta}}$ be defined by
(\ref{defhattheta}) and (\ref{multdeftheta}), respectively.
Then
\[
\hat{\theta}_{\mathbf{i}} =\theta_{\mathbf{i}}+ O_p
\bigl(n^{-1/2}\delta_n^{-|\mathbf{i}|} \bigr)+ O \bigl(
\delta_n^{\alpha-|\mathbf{i}|} \bigr) + O \bigl(\|\tilde{\bolds{\mu}}-
\bolds{\mu}\|^\alpha \delta_n^{-|\mathbf{i}|} \bigr),\qquad
\mathbf{i} \in \mathbb{I}.
\]
\end{lemma}

\begin{pf}
Using (\ref{defhattheta}) and Lemma~\ref{apprlemmult}, write
\begin{equation}
\label{diffmult} \hat{{\bolds{\theta}}}-{\bolds{\theta}}= \bigl(
\mathbf{Z}^T \mathbf{Z}\bigr)^{-1}\mathbf{Z}^T
\mathbf{Y} -{\bolds{\theta}} = \bigl(\mathbf{Z}^T \mathbf{Z}
\bigr)^{-1} \mathbf{Z}^T(\bolds{\eta}+ \bolds{\xi}).
\end{equation}
Since
$\mathrm{E}(\bolds{\xi})=\mathbf{0}$ and
$\operatorname{Cov} ((\mathbf{Z}^T \mathbf{Z})^{-1} \mathbf{Z}^T \bolds{\xi} )
= \sigma^2 (\mathbf{Z}^T \mathbf{Z})^{-1}$,  the order of the term
$(\mathbf{Z}^T \mathbf{Z})^{-1}\mathbf{Z}^T\bolds{\xi}$ is determined by
the diagonal entries of the matrix $(\mathbf{Z}^T \mathbf{Z})^{-1}$.
Hence, by Lemma~\ref{Zlemmult},  we have
\begin{equation}
\label{epsconvmult} \qquad\bigl(\mathbf{Z}^T
\mathbf{Z}\bigr)^{-1}\mathbf{Z}^T \bolds{\xi} =
\bigl(O_p \bigl(n^{-1/2} \bigr), O_p
\bigl(n^{-1/2}\delta_n^{-|\mathbf{i}_1|} \bigr), \ldots,
\\
O_p \bigl(n^{-1/2}\delta_n^{-|\mathbf{i}_q|}
\bigr) \bigr)^T.
\end{equation}

In view of (\ref{designmult}), $\mathbf{z}_k \in C(l\delta_n)$,  so
that
$|\mathbf{z}_k^{\mathbf{i}}| \le c \delta_n^{|\mathbf{i}|}$, $k=1,\ldots, n_2$,
$\mathbf{i}\in\mathbb{I}$.
Using this, (\ref{multeta}), $n_2\le n$
and Lemma~\ref{Zlemmult}, we obtain that
\begin{eqnarray}\label{etaconvmult}
\bigl(\mathbf{Z}^T \mathbf{Z}\bigr)^{-1}
\mathbf{Z}^T \bolds{\eta} & =& \pmatrix{ h_{00}\displaystyle \sum
_{k=1}^{n_2} \mathbf{z}_k^{\mathbf{i}_0}
\eta_k + \cdots + h_{0q} \displaystyle\sum
_{k=1}^{n_2} \mathbf{z}_k^{\mathbf{i}_q}
\eta_k \vspace*{2pt}
\cr
h_{10} \displaystyle\sum
_{k=1}^{n_2} \mathbf{z}_k^{\mathbf{i}_0}
\eta_k + \cdots + h_{1q} \displaystyle\sum_{i=k}^{n_2}
\mathbf{z}_k^{\mathbf{i}_q} \eta_k\vspace*{2pt}
\cr
\vdots
\vspace*{2pt}
\cr
h_{q0} \displaystyle\sum_{k=1}^{n_2}
\mathbf{z}_k^{\mathbf{i}_0} \eta_k + \cdots +
h_{qq} \displaystyle\sum_{k=1}^{n_2}
\mathbf{z}_k^{\mathbf{i}_q} \eta_k }
\nonumber
\\[-8pt]
\\[-8pt]
\nonumber
& =& \pmatrix{ O\bigl(
\delta_n^\alpha\bigr) + O \bigl(\|\tilde{\bolds{\mu}}-\bolds{
\mu}\|^\alpha \bigr) \vspace*{2pt}
\cr
O \bigl(\delta_n^{\alpha-|\mathbf{i}_1|}
\bigr) + O \bigl(\|\tilde{\bolds{\mu}}-\bolds{\mu}\|^\alpha
\delta_n^{-|\mathbf{i}_1|} \bigr)\vspace*{2pt}
\cr
\vdots \vspace*{2pt}
\cr
O \bigl(\delta_n^{\alpha-|\mathbf{i}_q|} \bigr) + O \bigl(\|\tilde{\bolds{
\mu}}-\bolds{\mu}\|^\alpha\delta_n^{-|\mathbf{i}_q|} \bigr) }.
\end{eqnarray}
Combining relations (\ref{diffmult}), (\ref{epsconvmult}) and (\ref{etaconvmult})
completes the proof of the lemma.
\end{pf}

Let $\mathbf{1}_j$, $j=1,\ldots, d$, be the $d$ standard unit vectors of $\mathbb{R}^d$,
that is,  $\mathbf{1}_j$ has $1$ at the $j$th coordinate and zeros
at other $d-1$ coordinates.
Notice that $\mathbb{I}_1 = \{\mathbf{1}_j\dvtx j=1,\ldots, d\}$.

\begin{lemma}
\label{nablalemmult}
Assume \textup{(A1)} and let  $f_{\bolds{\vartheta}}$, $\hat{{\bolds{\theta}}}$
and  ${\bolds{\theta}}$ be defined by
(\ref{multpolynomial}), (\ref{defhattheta})  and
(\ref{multdeftheta}),
respectively. If $\|\bolds{\mu}-\tilde{\bolds{\mu}}\|=o_p(\delta_n)$, then
\[
\nabla f_{\hat{{\bolds{\theta}}}}(\bolds{\mu}-\tilde{\bolds{\mu}}) -\nabla
f_{{\bolds{\theta}}}(\bolds{\mu}-\tilde{\bolds{\mu}}) =O_p
\bigl(n^{-1/2} \delta_n^{-1} \bigr) +
O_p \bigl(\delta_n^{\alpha-1} \bigr).
\]
\end{lemma}

\begin{pf}
The $j$th coordinate of the vector $\nabla f_{{\bolds{\theta}}}(\mathbf{x})$ is
\begin{eqnarray*}
\frac{\partial f_{{\bolds{\theta}}}(\mathbf{x})}{\partial x_j} &=& \sum_{\mathbf{i}\in\mathbb{I}} {\bolds{
\theta}}_{\mathbf{i}} \frac{\partial \mathbf{x}^{\mathbf{i}}}{\partial x_j} =\sum_{\mathbf{i}\in\mathbb{I}\dvtx i_j\ge 1}
{\bolds{\theta}}_{\mathbf{i}} \frac{\partial \mathbf{x}^{\mathbf{i}}}{\partial x_j} = \sum
_{\mathbf{i}\in\mathbb{I}\dvtx i_j\ge 1} i_j \theta_{\mathbf{i}}
\mathbf{x}^{\mathbf{i} - \mathbf{1}_j}
\\
&=& \theta_{\mathbf{1}_j} +\sum_{\mathbf{i}\in\mathbb{I}\dvtx i_j\ge 1, |\mathbf{i}|\ge 2}
i_j \theta_{\mathbf{i}} \mathbf{x}^{\mathbf{i} - \mathbf{1}_j},
\end{eqnarray*}
where $\mathbf{1}_j  \in \mathbb{I}_1$.
Then, for each $j=1,\ldots, d$,
\begin{eqnarray}\label{nablak}
&&\frac{\partial f_{\hat{{\bolds{\theta}}}}
(\bolds{\mu} -\tilde{\bolds{\mu}})}{\partial x_j} - \frac{\partial f_{{\bolds{\theta}}}(\bolds{\mu}
-\tilde{\bolds{\mu}})}{\partial x_j}
\nonumber
\\[-8pt]
\\[-8pt]
\nonumber
&&\qquad=  \hat{\theta}_{\mathbf{1}_j} -
\theta_{\mathbf{1}_j} + \sum_{\mathbf{i}\in\mathbb{I}\dvtx i_j\ge 1, |\mathbf{i}|\ge 2} i_j
(\hat{\theta}_{\mathbf{i}} -\theta_{\mathbf{i}}) (\bolds{\mu}-\tilde{\bolds{
\mu}})^{\mathbf{i} - \mathbf{1}_j}.
\end{eqnarray}

Now we bound the right-hand side of (\ref{nablak}).
Since $\mathbf{1}_j\in \mathbb{I}_1$, that is, $|\mathbf{1}_j|=1$ for $j=1,\ldots, d$, and
$\|\bolds{\mu}-\tilde{\bolds{\mu}}\|=o_p(\delta_n)$,
we obtain by Lemma~\ref{thetalemmult} that
\begin{eqnarray}\label{new11}
\hat{\theta}_{\mathbf{1}_j} -\theta_{\mathbf{1}_j} &=& O_p
\bigl(n^{-1/2}\delta_n^{-1} \bigr)+ O \bigl(
\delta_n^{\alpha-1} \bigr) + O \bigl(\|\tilde{\bolds{\mu}}-
\bolds{\mu}\|^\alpha \delta_n^{-1} \bigr)
\nonumber
\\[-8pt]
\\[-8pt]
\nonumber
&= & O_p \bigl(n^{-1/2}\delta_n^{-1}
\bigr)+ O_p \bigl(\delta_n^{\alpha-1} \bigr), \qquad j=1,
\ldots, d.
\end{eqnarray}
The same argument applies to each term of the sum in
the right-hand side of (\ref{nablak}):
for all $\mathbf{i}\in\mathbb{I}$ such that  $i_j\ge 1$ and
$|\mathbf{i}|\ge 2$
\begin{eqnarray}\label{new12}
&&\bigl|(\hat{\theta}_{\mathbf{i}} - \theta_{\mathbf{i}}) (\bolds{\mu}-\tilde{
\bolds{\mu}})^{\mathbf{i} - \mathbf{1}_j}\bigr |
\nonumber
\\
&&\qquad\le |\hat{\theta}_{\mathbf{i}} -\theta_{\mathbf{i}}| \|\bolds{\mu}-\tilde{
\bolds{\mu}}\|^{|\mathbf{i} - \mathbf{1}_j|}
\nonumber
\\
&&\qquad= \bigl[ O_p \bigl(n^{-1/2}\delta_n^{-|\mathbf{i}|}
\bigr)+ O \bigl(\delta_n^{\alpha-|\mathbf{i}|} \bigr) + O \bigl(\|\tilde{
\bolds{\mu}}-\bolds{\mu}\|^\alpha \delta_n^{-|\mathbf{i}|}
\bigr) \bigr] \|\bolds{\mu}-\tilde{\bolds{\mu}}\|^{|\mathbf{i}|-1}
\\
&&\qquad= o_p \bigl(n^{-1/2}\delta_n^{-1}
\bigr)+ o_p \bigl(\delta_n^{\alpha-1} \bigr) +
o_p \bigl(\|\tilde{\bolds{\mu}}-\bolds{\mu}\|^\alpha
\delta_n^{-1} \bigr)
\nonumber
\\
&&\qquad= o_p \bigl(n^{-1/2}\delta_n^{-1}\nonumber
\bigr)+ o_p \bigl(\delta_n^{\alpha-1} \bigr).
\end{eqnarray}
There are fixed number of terms in the sum from (\ref{nablak})
and the constant $i_j$ is at most  $r_\alpha$. Combining this with
(\ref{new11}) and (\ref{new12}), we see that the main term in (\ref{nablak}) is
$\hat{\theta}_{\mathbf{1}_j} -\theta_{\mathbf{1}_j}$ and therefore
\[
\bigl\|\nabla f_{\hat{{\bolds{\theta}}}}(\bolds{\mu}-\tilde{\bolds{\mu}}) -\nabla
f_{{\bolds{\theta}}}(\bolds{\mu}-\tilde{\bolds{\mu}})\bigr\| =O_p
\bigl(n^{-1/2} \delta_n^{-1} \bigr) +
O_p \bigl(\delta_n^{\alpha-1} \bigr).
\]
\upqed\end{pf}

For an $(s\times p)$-matrix $\mathbf{A}$, let
$\|\mathbf{A}\| = \sup_{\mathbf{x} \in \mathbb{R}^p\dvtx \|\mathbf{x}\| \le 1} \|\mathbf{A}\mathbf{x}\|$
be the operator norm for the rest of this section
and define the maximum norm $\|\mathbf{A}\|_{\max}=\max_{i,j} |a_{ij}|$, where
$a_{ij}$ are the entries of the matrix $\mathbf{A}$.
These norms are related by
\begin{equation}
\label{normsrelation} \|\mathbf{A}\|_{\max} \le \|
\mathbf{A}\| \le \sqrt{sp} \|\mathbf{A}\|_{\max}.
\end{equation}

\begin{lemma}
\label{Hessianlemmult}
Assume \textup{(A1)}, \textup{(A2)}, $\|\bolds{\mu}-\tilde{\bolds{\mu}}\|
= o_p(\delta_n)$ and $\sqrt{n}\delta_n^2\to\infty$.
For $\bolds{\mu}^*\in \mathbb{R}^d$ such that
$\|\bolds{\mu}^*\|=o_p(1)$ and for any fixed $\varepsilon \in (0,1)$,
let
\begin{equation}
\label{new13} B_n= \bigl\{\bigl\|Hf(\bolds{\mu}) -Hf_{\hat{{\bolds{\theta}}}}
\bigl(\bolds{\mu}^*\bigr)\bigr\| \le (1-\varepsilon) \bigl\| \bigl(Hf(\bolds{\mu})
\bigr)^{-1} \bigr\|^{-1} \bigr\}.
\end{equation}
Then $\P(B_n) \to 1$ as $n\to \infty$, on the event $B_n$,
$ (Hf_{\hat{{\bolds{\theta}}}}(\bolds{\mu}^*) )^{-1}$ exists and
\[
\bigl\| \bigl(Hf(\bolds{\mu}) \bigr)^{-1} - \bigl(Hf_{\hat{{\bolds{\theta}}}}\bigl(
\bolds{\mu}^*\bigr) \bigr)^{-1} \bigr\| = o_p(1).
\]
\end{lemma}

\begin{pf}
Clearly, by the smoothness of a polynomial,
\begin{equation}
\label{rel1Hessianlem} Hf_{{\bolds{\theta}}}(
\mathbf{z}) = Hf_{{\bolds{\theta}}}(\mathbf{0}) + O\bigl(\|\mathbf{z}\|\bigr)\qquad \mbox{as } \|
\mathbf{z}\| \to 0.
\end{equation}
We note that the elements of the matrix $Hf_{{\bolds{\theta}}}(\mathbf{0})$
[resp., $Hf_{\hat{\bolds{\theta}}}(\mathbf{0})$] are linear combinations of
${\theta}_{\mathbf{i}}$ (resp., $\hat{\theta}_{\mathbf{i}}$), $\mathbf{i} \in \mathbb{I}_2$.
From Lemma~\ref{thetalemmult} and the conditions
$\alpha>2$, $\|\bolds{\mu}-\tilde{\bolds{\mu}}\| = o_p(\delta_n)$
and $\sqrt{n}\delta_n^2\to\infty$, we obtain that
\[
\hat{\theta}_{\mathbf{i}}- \theta_{\mathbf{i}} = O_p
\bigl(n^{-1/2}\delta_n^{-2} \bigr)+ O \bigl(
\delta_n^{\alpha-2} \bigr) + O \bigl(\|\tilde{\bolds{\mu}}-
\bolds{\mu}\|^\alpha \delta_n^{-2} \bigr) =
o_p(1),\qquad  \mathbf{i} \in \mathbb{I}_2,
\]
where vector ${\bolds{\theta}}$ is defined by (\ref{multdeftheta}).
Therefore, entry-wise
\begin{equation}
\label{rel3Hessianlem} Hf_{\hat{{\bolds{\theta}}}}(
\mathbf{0}) =Hf_{\bolds{\theta}}(\mathbf{0}) +o_p(1).
\end{equation}
By (A1), (A2) and the definition (\ref{multdeftheta}) of ${\bolds{\theta}}$,
$Hf(\bolds{\mu}) = Hf_{{\bolds{\theta}}}(\bolds{\mu}-\tilde{\bolds{\mu}})$.
This and (\ref{rel1Hessianlem})
imply that entry-wise
\begin{equation}
\label{rel4Hessianlem} Hf(\bolds{\mu}) =
Hf_{{\bolds{\theta}}}(\bolds{\mu}-\tilde{\bolds{\mu}}) = Hf_{\bolds{\theta}}(
\mathbf{0})+O\bigl(\|\bolds{\mu}-\tilde{\bolds{\mu}}\|\bigr).
\end{equation}
Combining (\ref{rel1Hessianlem}), (\ref{rel3Hessianlem})
and (\ref{rel4Hessianlem}) leads to the following entry-wise relation:
\begin{eqnarray*}
Hf_{\hat{{\bolds{\theta}}}}\bigl(\bolds{\mu}^*\bigr)&=& Hf_{\hat{{\bolds{\theta}}}}(\mathbf{0})+O
\bigl(\|\bolds{\mu}^*\|\bigr)
\\
&=& Hf_{\bolds{\theta}}(\mathbf{0}) +o_p(1)+O\bigl(\|\bolds{\mu}^*\|
\bigr)
\\
&=& Hf(\bolds{\mu}) +O\bigl(\|\bolds{\mu}-\tilde{\bolds{\mu}}\|\bigr)+o_p(1)+O
\bigl(\|\bolds{\mu}^*\|\bigr)
\\
&=& Hf(\bolds{\mu}) +o_p(1).
\end{eqnarray*}
Then $\|Hf_{\hat{{\bolds{\theta}}}}(\bolds{\mu}^*)-Hf(\bolds{\mu})\|_{\max} = o_p(1)$
and hence, by (\ref{normsrelation}),
\begin{equation}
\label{rel4aHessianlem}\bigl \|Hf(\bolds{\mu})
-Hf_{\hat{{\bolds{\theta}}}}\bigl(\bolds{\mu}^*\bigr)\bigr\| = o_p(1).
\end{equation}

Next, since (A1) and (A2) imply (\ref{newA3}),
$\lambda_{\min}(Hf(\bolds{\mu}))
\le \cdots \le \lambda_{\max}(Hf(\bolds{\mu})) \le -\lambda_0 <0$.
Hence,
$ \| (Hf(\bolds{\mu}) )^{-1} \|=
- (\lambda_{\max}(Hf(\bolds{\mu})) )^{-1}\le
\lambda_0^{-1}$,  or
\begin{equation}
\label{rel5Hessianlem} \lambda_0 \le \bigl\|
\bigl(Hf(\bolds{\mu}) \bigr)^{-1} \bigr\|^{-1}.
\end{equation}

Define the event $C_n= \{\|Hf(\bolds{\mu}) -
Hf_{\hat{{\bolds{\theta}}}}(\bolds{\mu}^*)\| \le (1-\varepsilon) \lambda_0  \}$.
Using (\ref{rel5Hessianlem}) and Lemma~\ref{lem2appendix},
we obtain that
\[
C_n \subset B_n\subset \bigl\{ \bigl(Hf_{\hat{{\bolds{\theta}}}}
\bigl(\bolds{\mu}^*\bigr) \bigr)^{-1} \mbox{ exists} \bigr\}.
\]
In view  of (\ref{rel4aHessianlem}), $\P(C_n) \to 1$ and hence $\P(B_n) \to 1$.
Finally, by applying (\ref{rel4aHessianlem}),  (\ref{rel5Hessianlem})
and Lemma~\ref{lem2appendix} again, we get that
on the event $B_n$
\begin{eqnarray*}
\bigl\| \bigl(Hf(\bolds{\mu}) \bigr)^{-1} - \bigl(Hf_{\hat{{\bolds{\theta}}}}\bigl(
\bolds{\mu}^*\bigr) \bigr)^{-1}\bigr \| &\le& \varepsilon^{-1} \bigl\|
\bigl(Hf(\bolds{\mu}) \bigr)^{-1}\bigr \|^2\bigl \|Hf(\bolds{\mu})
-Hf_{\hat{{\bolds{\theta}}}}\bigl(\bolds{\mu}^*\bigr)\bigr\|
\\
&\le& \varepsilon^{-1} \lambda_0^{-2}\bigl \|Hf(\bolds{
\mu}) -Hf_{\hat{{\bolds{\theta}}}}\bigl(\bolds{\mu}^*\bigr)\bigr\| = o_p(1).
\end{eqnarray*}
\upqed\end{pf}

\begin{remark}
Lemma~\ref{Hessianlemmult} would still hold if we only
assumed that $\|\bolds{\mu}-\tilde{\bolds{\mu}}\|
= O_p(\delta_n)$ instead of  $\|\bolds{\mu}-\tilde{\bolds{\mu}}\|
= o_p(\delta_n)$.
\end{remark}

\begin{lemma}
\label{Alemmult}
Assume \textup{(A1)}, \textup{(A2)}, $\|\bolds{\mu}-\tilde{\bolds{\mu}}\| = o_p(\delta_n)$,
$\sqrt{n}\delta_n^2\to\infty$
and let $A_n = \{\mathring{\bolds{\mu}}\in C(2l\delta_n/3) \}$,
where the estimator $\mathring{\bolds{\mu}}$ is defined by (\ref{muestim2}).
Then  $\P(A_n)\to 1$ as $n\to \infty$.
\end{lemma}

\begin{pf}
Bound $\P({A}_n^c)$ by
\begin{equation}
\label{rel1lem7} \P \bigl(\mathring{\bolds{\mu}}\notin C(2l
\delta_n/3), \bolds{\mu}- \tilde{\bolds{\mu}} \in C(l
\delta_n/3) \bigr) +\P \bigl(\bolds{\mu}- \tilde{\bolds{\mu}} \notin
C(l\delta_n/3) \bigr).
\end{equation}
The second term converges to zero
by the condition  $\|\bolds{\mu}-\tilde{\bolds{\mu}}\| = o_p(\delta_n)$.

For a symmetric matrix $\mathbf{M}$ and any $\mathbf{x}\in\mathbb{R}^d$,
$\lambda_{\min}(\mathbf{M}) \|\mathbf{x}\|^2 \le \mathbf{x}^T \mathbf{M} \mathbf{x}
\le   \lambda_{\max}(\mathbf{M}) \|\mathbf{x}\|^2$.
Recall that $\nabla f(\bolds{\mu})=\mathbf{0}$ and (\ref{newA3})
follow from (A1) and (A2).
Then, for $\bolds{\mu}\in C(\tilde{\bolds{\mu}},l\delta_n/3)$
and $\mathbf{x}\in C(\tilde{\bolds{\mu}},l\delta_n)\setminus C(\tilde{\bolds{\mu}},2l\delta_n/3)$,
by using Taylor's expansion, $\nabla f(\bolds{\mu})=\mathbf{0}$ and (\ref{newA3}), we have
\begin{eqnarray}\label{lem1rel1}
f(\mathbf{x}) &= &f(\bolds{\mu}) + \frac12 (\mathbf{x}-\bolds{\mu})^T
Hf\bigl(\bolds{\mu}^*\bigr) (\mathbf{x}-\bolds{\mu})
\nonumber
\\[-8pt]
\\[-8pt]
\nonumber
 &\le& f(\bolds{\mu}) - \frac{\lambda_0}{2} \|
\mathbf{x}-\bolds{\mu}\|^2 \le f(\bolds{\mu}) - c
\delta_n^2
\end{eqnarray}
for some positive constant $c=c(\lambda_0,l)$
and sufficiently large $n$ such that $\|\bolds{\mu}^* - \bolds{\mu}\| \le \kappa$, with
$\kappa>0$ from (\ref{newA3}).

Next, by using (\ref{multHolder}), (\ref{multdeftheta}) and the
$c_r$-inequality,
\[
f_{{\bolds{\theta}}}(\mathbf{z}) = P_{f,\bolds{\mu}}(\mathbf{z}+\bolds{\tilde{
\mu}})= f(\mathbf{z}+\tilde{\bolds{\mu}}) +O\bigl(\|\mathbf{z}\|^\alpha
\bigr)+ O\bigl(\|\tilde{\bolds{\mu}}-\bolds{\mu}\|^\alpha\bigr).
\]
Now we combine this with Lemma~\ref{thetalemmult}
and the conditions $\alpha>2$, $\|\bolds{\mu}-\tilde{\bolds{\mu}}\| = o_p(\delta_n)$
and $\sqrt{n}\delta_n^2\to\infty$
to obtain that, uniformly in  $\mathbf{z} \in C(l\delta_n)$,
\begin{eqnarray}\label{lem1rel2}
f_{\hat{{\bolds{\theta}}}}(\mathbf{z}) &=& f_{\bolds{\theta}}(\mathbf{z}) +O_p
\bigl( n^{-1/2} \bigr)+O\bigl(\delta_n^\alpha\bigr)
\nonumber
\\
&=& f(\mathbf{z}+\tilde{\bolds{\mu}}) +O\bigl(\|\mathbf{z}\|^\alpha
\bigr)+ O \bigl(\|\tilde{\bolds{\mu}}-\bolds{\mu}\|^\alpha \bigr)+
O_p \bigl( n^{-1/2} \bigr) +O\bigl(\delta_n^\alpha
\bigr)
\\
&=&  f(\mathbf{z}+\tilde{\bolds{\mu}}) +
o_p\bigl(\delta_n^2\bigr).\nonumber
\end{eqnarray}

Recall that $\mathring{\bolds{\mu}} \in C(l\delta_n)$ by the
definition (\ref{muestim2}). By (\ref{lem1rel1}) and
(\ref{lem1rel2}), we see that the event
\begin{eqnarray*}
&&\bigl\{\mathring{\bolds{\mu}} \notin C(2l\delta_n/3), \bolds{\mu}
 - \tilde{\bolds{\mu}} \in C(l\delta_n/3)\bigr\}
\\
&&\qquad= \bigl\{\mathring{\bolds{\mu}}+\tilde{\bolds{\mu}} \in C(\tilde{\bolds{\mu}},l
\delta_n)\setminus C(\tilde{\bolds{\mu}},2l\delta_n/3),
\bolds{\mu} \in C(\tilde{\bolds{\mu}},l\delta_n/3)\bigr\}
\end{eqnarray*}
implies the event
\begin{eqnarray*}
f(\bolds{\mu}) - c \delta_n^2 & \ge& f(\mathring{\bolds{
\mu}}+\tilde{\bolds{\mu}}) = f_{\hat{{\bolds{\theta}}}}(\mathring{\bolds{
\mu}})+o_p\bigl(\delta_n^2\bigr)
\\
& \ge& f_{\hat{{\bolds{\theta}}}}(\bolds{\mu}-\tilde{\bolds{\mu}})+o_p\bigl(
\delta_n^2\bigr)= f(\bolds{\mu}) +o_p\bigl(
\delta_n^2\bigr),
\end{eqnarray*}
leading to
\[
\P \bigl(\mathring{\bolds{\mu}}\notin C(2l\delta_n/3), \bolds{
\mu}- \tilde{\bolds{\mu}} \in C(l\delta_n/3) \bigr) \le \P\bigl(c
\delta_n^2 \le o_p\bigl(\delta_n^2
\bigr)\bigr) \to 0
\]
as $n\to \infty$. Combined with (\ref{rel1lem7}), this completes the
proof of the lemma.
\end{pf}

\begin{pf*}{Proof of Theorem~\ref{maintheorem}}
By (A1) and (A2), $\nabla f(\bolds{\mu}) =\mathbf{0}$.
According to the definition (\ref{multdeftheta}) of the polynomial
$f_{\bolds{\theta}}$,
\begin{equation}
\label{rel1theorem1} \mathbf{0} = \nabla f(\bolds{\mu})= \nabla
P_{f,\bolds{\mu}}(\bolds{\mu}) = \nabla f_{\bolds{\theta}}(\bolds{\mu}-\tilde{
\bolds{\mu}}).
\end{equation}

By (\ref{muestim2}), $\max_{\mathbf{z}\in C(l\delta_n)}
f_{\hat{{\bolds{\theta}}}}(\mathbf{z})=f_{\hat{{\bolds{\theta}}}}(\mathring{\bolds{\mu}})$.
If this maximum is not attained on the boundary of $C(l\delta_n)$,
then $\nabla f_{\hat{{\bolds{\theta}}}}(\mathring{\bolds{\mu}})$ must be zero.
Hence we have that on the event $A_n=\{\mathring{\bolds{\mu}}\in C(2l\delta_n/3)\}$
\begin{equation}
\label{rel2theorem1} \mathbf{0}=\nabla f_{\hat{{\bolds{\theta}}}}(
\mathring{\bolds{\mu}}) = \nabla f_{\hat{{\bolds{\theta}}}}(\bolds{\mu}-\tilde{\bolds{
\mu}}) + Hf_{\hat{{\bolds{\theta}}}}\bigl(\bolds{\mu}^*\bigr) \bigl(\mathring{\bolds{\mu}}-
(\bolds{\mu}-\tilde{\bolds{\mu}})\bigr),
\end{equation}
where $\bolds{\mu}^*=(\mu^*_1,\ldots, \mu^*_d)
=\lambda \mathring{\bolds{\mu}}+(1-\lambda)(\bolds{\mu}-\tilde{\bolds{\mu}})$
for some $\lambda \in[0,1]$.
Thus  $\|\bolds{\mu}^*\| = O(\|\mathring{\bolds{\mu}}\|) +
O(\|\bolds{\mu}-\tilde{\bolds{\mu}}\|)= O_p(\delta_n)=o_p(1)$.

By Lemma~\ref{Hessianlemmult},
$ (Hf_{\hat{{\bolds{\theta}}}}(\bolds{\mu}^*) )^{-1}$ exists on the event
$B_n$ defined by (\ref{new13}).
Relations~(\ref{rel1theorem1}) and (\ref{rel2theorem1}) imply
that  on the event $A_n\cap B_n$
\begin{eqnarray} \label{rel3theorem1}
\hat{\bolds{\mu}}- \bolds{\mu} &=& - \bigl(Hf_{\hat{{\bolds{\theta}}}}\bigl(\bolds{\mu}^*
\bigr) \bigr)^{-1} \nabla f_{\hat{{\bolds{\theta}}}}(\bolds{\mu}-\tilde{\bolds{
\mu}})
\nonumber
\\
&=& - \bigl(Hf_{\hat{{\bolds{\theta}}}}\bigl(\bolds{\mu}^*\bigr) \bigr)^{-1}
\bigl(\nabla f_{\hat{{\bolds{\theta}}}}(\bolds{\mu}-\tilde{\bolds{\mu}}) -\nabla
f_{{\bolds{\theta}}}(\bolds{\mu}-\tilde{\bolds{\mu}}) \bigr)
\\
&=& - \bigl(Hf(\bolds{\mu})
\bigr)^{-1} \bigl(\nabla f_{\hat{{\bolds{\theta}}}}(\bolds{\mu}-\tilde{\bolds{
\mu}}) -\nabla f_{{\bolds{\theta}}}(\bolds{\mu}-\tilde{\bolds{\mu}}) \bigr) +
r_n,\nonumber
\end{eqnarray}
where
$r_n =  [ (Hf(\bolds{\mu}) )^{-1} -
(Hf_{\hat{{\bolds{\theta}}}}(\bolds{\mu}^*) )^{-1}  ]
(\nabla f_{\hat{{\bolds{\theta}}}}(\bolds{\mu}-\tilde{\bolds{\mu}})
-\nabla f_{{\bolds{\theta}}}(\bolds{\mu}-\tilde{\bolds{\mu}}) )$ is the remainder term.

By Lemma~\ref{nablalemmult} and (\ref{rel5Hessianlem}),
we bound the norm of the first term on the right-hand side of
(\ref{rel3theorem1}) as
\begin{eqnarray*}
&&\bigl\| \bigl(Hf(\bolds{\mu}) \bigr)^{-1} \bigl(\nabla f_{\hat{{\bolds{\theta}}}}(
\bolds{\mu}-\tilde{\bolds{\mu}})  -\nabla f_{{\bolds{\theta}}}(\bolds{\mu}-\tilde{
\bolds{\mu}}) \bigr) \bigr\|
\\
&&\qquad \le \lambda_0^{-1} \bigl\|\nabla f_{\hat{{\bolds{\theta}}}}(\bolds{
\mu}-\tilde{\bolds{\mu}}) -\nabla f_{{\bolds{\theta}}}(\bolds{\mu}-\tilde{\bolds{
\mu}})\bigr\|=O_p(\gamma_n),
\end{eqnarray*}
where $\gamma_n =n^{-1/2}\delta_n^{-1}+\delta_n^{\alpha-1}$.
Therefore $\|r_n\| =o_p(1) O_p(\gamma_n)=o_p(\gamma_n)$
on the event $B_n$ by Lemmas~\ref{nablalemmult}
and~\ref{Hessianlemmult}. Consequently on the event $A_n\cap B_n$,
we have
\begin{equation}
\label{new14} \|\hat{\bolds{\mu}}- \bolds{\mu} \|= O_p
\bigl(n^{-1/2}\delta_n^{-1}+\delta_n^{\alpha-1}
\bigr)= O_p(\gamma_n).
\end{equation}

For any constant $\rho>0$,
\[
\P \bigl(\|\hat{\bolds{\mu}} -\bolds{\mu}\| >\rho \gamma_n \bigr) \le \P
\bigl(\bigl\{\|\hat{\bolds{\mu}} -\bolds{\mu}\| >\rho \gamma_n\bigr\}\cap
A_n\cap B_n \bigr) + \P\bigl({A}_n^c
\bigr) +\P\bigl({B}_n^c\bigr).
\]
The first term on the right-hand side can  be made arbitrarily small
by choosing~$\rho$ sufficiently large in view of
(\ref{new14}), uniformly in $n$, while the other two terms converge to zero by
Lemmas~\ref{Hessianlemmult} and~\ref{Alemmult}.
This proves (\ref{new5}).

It remains to prove (\ref{new6}).
From (\ref{multdeftheta}) it follows that
\[
M=f(\bolds{\mu}) = f_{\bolds{\theta}}(\bolds{\mu}-\tilde{\bolds{\mu}}) = \sum
_{\mathbf{i}\in\mathbb{I}} \theta_{\mathbf{i}}(\bolds{\mu}-\tilde{
\bolds{\mu}})^{\mathbf{i}},
\]
so that, according to (\ref{muestim2}), $\hat{M} - M$ can be written as
\begin{eqnarray} \label{rel5theorem1}
&&f_{\hat{{\bolds{\theta}}}}(\mathring{\bolds{\mu}})  - f_{\bolds{\theta}}(\bolds{\mu}-
\tilde{\bolds{\mu}})\nonumber\\
 &&\qquad= \sum_{\mathbf{i}\in\mathbb{I}} \bigl[\hat{
\theta}_{\mathbf{i}} \mathring{\bolds{\mu}}^{\mathbf{i}} -
\theta_{\mathbf{i}} (\bolds{\mu}-\tilde{\bolds{\mu}})^{\mathbf{i}} \bigr]
\\
&&\qquad= \hat{\theta}_{\mathbf{i}_0} -
\theta_{\mathbf{i}_0} +\sum_{\mathbf{i}\in\mathbb{I}, |\mathbf{i}|\ge 1} (\hat{
\theta}_{\mathbf{i}}-\theta_{\mathbf{i}}) \mathring{\bolds{
\mu}}^{\mathbf{i}} +\sum_{\mathbf{i}\in\mathbb{I}, |\mathbf{i}|\ge 1}
\theta_{\mathbf{i}} \bigl[\mathring{\bolds{\mu}}^{\mathbf{i}} -(\bolds{\mu}-
\tilde{\bolds{\mu}})^{\mathbf{i}} \bigr].\nonumber
\end{eqnarray}
By Lemma~\ref{thetalemmult}, the first term in (\ref{rel5theorem1}) is
\begin{equation}
\label{rel6theorem1} \hat{\theta}_{\mathbf{i}_0} -
\theta_{\mathbf{i}_0} = O_p\bigl(n^{-1/2}\bigr) +
O_p\bigl(\delta_n^\alpha\bigr).
\end{equation}
From (\ref{muestim2}), (\ref{new5}) and the conditions
$\sqrt{n} \delta_n^2 \to \infty$,
$\alpha>2$, $\|\bolds{\mu} - \tilde{\bolds{\mu}}\|=o_p(\delta_n)$, it follows that
\begin{eqnarray}\label{rel6atheorem1}
\|\mathring{\bolds{\mu}} \| &=& \|\hat{\bolds{\mu}}-\tilde{\bolds{\mu}}\| \le \|
\hat{\bolds{\mu}}-\bolds{\mu}\| + \|\bolds{\mu} - \tilde{\bolds{\mu}}\|
\nonumber
\\[-8pt]
\\[-8pt]
\nonumber
&=&  O_p \bigl(n^{-1/2}
\delta_n^{-1} \bigr) +O_p \bigl(
\delta_n^{\alpha-1} \bigr) +O\bigl(\|\bolds{\mu} - \tilde{\bolds{\mu}}
\|\bigr) = o_p(\delta_n).
\end{eqnarray}

Using  (\ref{rel6atheorem1}) and Lemma~\ref{thetalemmult},
each term in the second sum of (\ref{rel5theorem1})
\[
\bigl|(\hat{\theta}_{\mathbf{i}}-\theta_{\mathbf{i}}) \mathring{\bolds{
\mu}}^{\mathbf{i}}\bigr| \le |\hat{\theta}_{\mathbf{i}}-\theta_{\mathbf{i}}| \|
\mathring{\bolds{\mu}}\|^{|\mathbf{i}|} =o_p\bigl(n^{-1/2}
\bigr) + o_p\bigl(\delta_n^\alpha\bigr),
\]
so that, as there are a fixed number of terms in the sum,
\begin{equation}
\label{rel7theorem1} \sum_{\mathbf{i}\in\mathbb{I}, |\mathbf{i}|\ge 1} (
\hat{\theta}_{\mathbf{i}}-\theta_{\mathbf{i}})\mathring{\bolds{
\mu}}^{\mathbf{i}} =o_p\bigl(n^{-1/2}\bigr) +
o_p\bigl(\delta_n^\alpha\bigr).
\end{equation}

Now consider the third sum in (\ref{rel5theorem1}).
Combining Lemma~\ref{lem0appendix} with  (\ref{new5}),
(\ref{rel6atheorem1}) and the condition
$\|\bolds{\mu} - \tilde{\bolds{\mu}}\|=o_p(\delta_n)$, we obtain that
for any  $\mathbf{i} \in \mathbb{I}$, $|\mathbf{i}|\ge 1$,
\begin{eqnarray} \label{rel8theorem1}
\bigl|\mathring{\bolds{\mu}}^{\mathbf{i}} -(\bolds{\mu}-\tilde{\bolds{
\mu}})^{\mathbf{i}}\bigr | &\le& \bigl\|\mathring{\bolds{\mu}}-(\bolds{\mu}-\tilde{\bolds{
\mu}})\bigr\| \sum_{k=1}^{|\mathbf{i}|}\|\mathring{\bolds{
\mu}}\|^{|\mathbf{i}| -k} \|\bolds{\mu}-\tilde{\bolds{\mu}}\|^{k-1}
\nonumber
\\
&=& \|\hat{\bolds{\mu}}-\bolds{\mu}\| \sum_{k=1}^{|\mathbf{i}|}
\|\mathring{\bolds{\mu}}\|^{|\mathbf{i}| -k} \|\bolds{\mu}-\tilde{\bolds{\mu}}
\|^{k-1}
\\
&=& o_p \bigl(n^{-1/2}
\delta_n^{|\mathbf{i}|-2} \bigr) + o_p \bigl(
\delta_n^{\alpha+|\mathbf{i}|-2} \bigr).\nonumber
\end{eqnarray}

Since $\mathbf{D}^{\mathbf{i}} P_{f,\bolds{\mu}}(\mathbf{x})$, $\mathbf{i}\in\mathbb{I}$,
are continuous, they are bounded over the compact set~$D$,
so that $\theta_{\mathbf{i}} = O_p(1)$, $\mathbf{i}\in\mathbb{I}$,
in view of  (\ref{new9}). Because of this  and (\ref{rel8theorem1}),
\begin{equation}
\label{rel9theorem1} \sum_{\mathbf{i}\in\mathbb{I}, |\mathbf{i}|\ge 2}
\theta_{\mathbf{i}} \bigl[\mathring{\bolds{\mu}}^{\mathbf{i}} -(\bolds{\mu}-
\tilde{\bolds{\mu}})^{\mathbf{i}} \bigr] =o_p\bigl(n^{-1/2}
\bigr) +o_p\bigl(\delta_n^\alpha\bigr).
\end{equation}

It remains to handle separately the terms in the third sum of
(\ref{rel5theorem1}) over $\mathbf{i} \in \mathbb{I}_1$, that is,
$\mathbf{i} \in \mathbb{I}$ such that $|\mathbf{i}|=1$.
Due to (\ref{new9}) and the condition
$\|\bolds{\mu} - \tilde{\bolds{\mu}}\|=o_p(\delta_n)$,
\begin{equation}
\label{new15} \theta_{\mathbf{i}} =\mathbf{D}^{\mathbf{i}}
P_{f,\bolds{\mu}}(\tilde{\bolds{\mu}}) = O\bigl(\|\bolds{\mu} - \tilde{\bolds{\mu}}
\|\bigr)= o_p(\delta_n), \qquad \mathbf{i} \in \mathbb{I}, |
\mathbf{i}|=1.
\end{equation}
Then (\ref{rel8theorem1}) and (\ref{new15}) imply that
\[
\sum_{\mathbf{i}\in\mathbb{I}, |\mathbf{i}|=1} \theta_{\mathbf{i}} \bigl[\mathring{
\bolds{\mu}}^{\mathbf{i}} -(\bolds{\mu}-\tilde{\bolds{\mu}})^{\mathbf{i}}
\bigr] =o_p\bigl(n^{-1/2}\bigr) +o_p\bigl(
\delta_n^\alpha\bigr).
\]

Finally, combining the last display with (\ref{rel5theorem1}),
(\ref{rel6theorem1}), (\ref{rel7theorem1}) and (\ref{rel9theorem1})
completes the proof of (\ref{new6}).
\end{pf*}

\begin{remark}
\label{remlast}
The above argument for estimating the parameter $M=f(\bolds{\mu})$ can be refined
for the problem of estimating any mixed derivative
$\mathbf{D}^{\mathbf{i}} f(\bolds{\mu})$, for $\mathbf{i} \in \mathbb{I}$, $|\mathbf{i}|\ge 2$.
One can  take the estimator
$\mathbf{D}^{\mathbf{i}} f_{\hat{{\bolds{\theta}}}}(\mathring{\bolds{\mu}})$
and establish in a similar way that
\[
\mathbf{D}^{\mathbf{i}} f_{\hat{{\bolds{\theta}}}}(\mathring{\bolds{\mu}})-
\mathbf{D}^{\mathbf{i}} f(\bolds{\mu}) = O_p \bigl(n^{-1/2}
\delta_n^{-|\mathbf{i}|} \bigr) + O_p \bigl(
\delta_n^{\alpha- |\mathbf{i}|} \bigr),\qquad  \mathbf{i} \in \mathbb{I}, |
\mathbf{i}|\ge 2.
\]
\end{remark}

\begin{appendix}\label{app}
\section*{Appendix}
\begin{lemma}
\label{lem0appendix}
For any $\mathbf{x},\mathbf{y}\in\mathbb{R}^d$ and
any $\mathbf{i} \in \mathbb{N}^d$ such that $|\mathbf{i}|\ge 1$,
\[
\bigl|\mathbf{x}^{\mathbf{i}}-\mathbf{y}^{\mathbf{i}}\bigr| \le \|\mathbf{x}-
\mathbf{y}\| \sum_{k=1}^{|\mathbf{i}|} \|\mathbf{x}
\|^{|\mathbf{i}| -k}\|\mathbf{y}\|^{k-1}.
\]
\end{lemma}

\begin{pf}
We prove the lemma by induction in dimension. For $d=1$,
\[
x^i - y^i =(x-y) \sum_{k=1}^i
x^{i-k} y^{k-1}
\]
and the statement follows.

Now we handle the inductive step.
Suppose the statement is true for all dimensions $k=1,\ldots, d-1$.
We want to show that it also holds for the dimension $d$. Without loss
of generality assume that  $i_1>0$. Recall the notation $\mathbf{x}_{-1}
=(x_2,\ldots, x_d)$, $\mathbf{i}_{-1} =(i_2,\ldots,i_d)$ that we used in
Lemma~\ref{indeplemmult}. We have
\begin{eqnarray*}
\mathbf{x}^{\mathbf{i}}-\mathbf{y}^{\mathbf{i}} &= &\mathbf{x}^{\mathbf{i}}
- x_1^{i_1} \mathbf{y}_{-1}^{\mathbf{i}_{-1}}+
x_1^{i_1} \mathbf{y}_{-1}^{\mathbf{i}_{-1}} -
\mathbf{y}^{\mathbf{i}}
\\
& = &x_1^{i_1}\bigl(\mathbf{x}_{-1}^{\mathbf{i}_{-1}}-
\mathbf{y}_{-1}^{\mathbf{i}_{-1}}\bigr) + \bigl(x_1^{i_1}-y_1^{i_1}
\bigr) \mathbf{y}_{-1}^{\mathbf{i}_{-1}}.
\end{eqnarray*}

Obviously, $|x_1|\le \|\mathbf{x}\|$,  $\|\mathbf{x}_{-1}\| \le \|\mathbf{x}\|$
and $|\mathbf{i}|=|\mathbf{i}_{-1}| +i_1$.
Using these relations and the assumption of the inductive step,
we  obtain that
\begin{eqnarray*}
\bigl|x_1^{i_1}\bigl(\mathbf{x}_{-1}^{\mathbf{i}_{-1}}-
\mathbf{y}_{-1}^{\mathbf{i}_{-1}}\bigr)\bigr| & \le& |x_1|^{i_1}
\|\mathbf{x}_{-1}-\mathbf{y}_{-1}\| \sum
_{k=1}^{|\mathbf{i}_{-1}|} \|\mathbf{x}_{-1}
\|^{|\mathbf{i}_{-1}| -k}\|\mathbf{y}_{-1}\|^{k-1}
\\
&\le& \|\mathbf{x}-\mathbf{y}\|\sum_{k=1}^{|\mathbf{i}|-i_1}
\|\mathbf{x}\|^{|\mathbf{i}| -k}\|\mathbf{y}\|^{k-1}
\end{eqnarray*}
and
\begin{eqnarray*}
\bigl|\bigl(x_1^{i_1}-y_1^{i_1}\bigr)
\mathbf{y}_{-1}^{\mathbf{i}_{-1}}\bigr| &\le& \|\mathbf{y}_{-1}
\|^{\mathbf{i}_{-1}} |x_1-y_1|\sum
_{k=1}^{i_1} |x_1|^{i_1-k}
|y_1|^{k-1}
\\
&\le& \|\mathbf{x}-\mathbf{y}\| \sum_{k=|\mathbf{i}|-i_1+1}^{|\mathbf{i}|}
\|\mathbf{x}\|^{|\mathbf{i}| -k}\|\mathbf{y}\|^{k-1}.
\end{eqnarray*}
Combining the last three relations, we obtain the desired result.\vadjust{\goodbreak}
\end{pf}

Below, we consider $s\times s$ matrices and let $\mathbf{I}$ denote
the identity matrix of order $s$. Let $\|\mathbf{A}\|$ be some norm on the space of
$s\times s$ matrices satisfying the multiplicative property
$\|\mathbf{A}\mathbf{B}\| \le \|\mathbf{A}\| \|\mathbf{B}\|$. For example, the
operator norm satisfies this property.

\begin{lemma}[(Banach's lemma)]
\label{Banachlem}
Let $\mathbf{M}$ be a matrix with $\|\mathbf{M}\| <1$.
Then $\mathbf{I}-\mathbf{M}$ is invertible,
$(\mathbf{I}-\mathbf{M})^{-1} = \mathbf{I}+\mathbf{M}+\mathbf{M}^2
+\cdots$ and $\|(\mathbf{I}-\mathbf{M})^{-1}\| \le (1-\|\mathbf{M}\|)^{-1}$.
\end{lemma}

The proof of Banach's lemma can be found in many textbooks on
functional analysis. The next two lemmas are essentially adopted from
Facer and M\"uller (\citeyear{FacerMuller2003}) with some modifications.

\begin{lemma}
\label{lem1appendix}
Let $\mathbf{V}$ be invertible and $\mathbf{W}$ be such that
$\|\mathbf{W}\| < \|\mathbf{V}^{-1}\|^{-1}$. Then $\mathbf{V}+\mathbf{W}$
is invertible and
\[
\bigl(\|\mathbf{V}\|+\|\mathbf{W}\|\bigr)^{-1} \le \bigl\|(\mathbf{V}+
\mathbf{W})^{-1}\bigr \| \le \frac{\|\mathbf{V}^{-1}\|}{1 - \|\mathbf{V}^{-1}\mathbf{W}\|}.
\]
\end{lemma}

\begin{pf}
Since $\|\mathbf{V}^{-1}\mathbf{W}\| <1$ due to the condition
$\|\mathbf{W}\| < \|\mathbf{V}^{-1}\|^{-1}$, the matrix
$(\mathbf{I}+\mathbf{V}^{-1}\mathbf{W})$ is invertible and
$\|(\mathbf{I}+\mathbf{V}^{-1}\mathbf{W})^{-1}\|\le (1 - \|\mathbf{V}^{-1}\mathbf{W}\|)^{-1}$
by Banach's lemma. Therefore,
$\mathbf{V}+\mathbf{W}=\mathbf{V}(\mathbf{I}+\mathbf{V}^{-1}\mathbf{W})$
is also invertible and
\begin{eqnarray*}
\bigl\|(\mathbf{V}+\mathbf{W})^{-1}\bigr\| &=& \bigl\|\bigl(\mathbf{I}+
\mathbf{V}^{-1}\mathbf{W}\bigr)^{-1}\mathbf{V}^{-1}
\bigr\|
\\
& \le& \bigl\|\mathbf{V}^{-1}\bigr\| \bigl\|\bigl(\mathbf{I}+\mathbf{V}^{-1}
\mathbf{W}\bigr)^{-1}\bigr\| \le \frac{\|\mathbf{V}^{-1}\|}{1 - \|\mathbf{V}^{-1}\mathbf{W}\|}.
\end{eqnarray*}

Now, using $\|\mathbf{V}+\mathbf{W}\| \le \|\mathbf{V}\| +\|\mathbf{W}\|$
and the invertibility of $\mathbf{V}+\mathbf{W}$,
we obtain $\|(\mathbf{V}+\mathbf{W})^{-1} \| \ge \|\mathbf{V}+\mathbf{W}\|^{-1}
\ge (\|\mathbf{V}\| +\|\mathbf{W}\|)^{-1}$.
\end{pf}

\begin{lemma}
\label{lem2appendix}
Let $\mathbf{A}$ be invertible and $\mathbf{B}$ be such that
$\| \mathbf{A}-\mathbf{B}\| \le (1-\varepsilon) \|\mathbf{A}^{-1}\|^{-1}$
for some $\varepsilon \in (0,1]$. Then $\mathbf{B}$ is invertible and
\[
\bigl\|\mathbf{B}^{-1}-\mathbf{A}^{-1}\bigr\| \le
\varepsilon^{-1} \bigl\|\mathbf{A}^{-1}\bigr\|^2 \| \mathbf{A}-
\mathbf{B}\|.
\]
\end{lemma}

\begin{pf}
Write $\mathbf{B}=\mathbf{A}+(\mathbf{B}-\mathbf{A})$ and apply Lemma~\ref{lem1appendix}
with $\mathbf{V}=\mathbf{A}$ and $\mathbf{W}=\mathbf{B}-\mathbf{A}$ to
conclude that $\mathbf{B}$ is invertible  and,
as $\|\mathbf{A}^{-1}(\mathbf{B}-\mathbf{A})\| \le 1-\varepsilon$ by
the condition of the lemma,
\[
\bigl\|\mathbf{B}^{-1}\bigr\| \le \frac{\|\mathbf{A}^{-1}\|}{1 -\|
\mathbf{A}^{-1}(\mathbf{B}-\mathbf{A})\|} \le \frac{\|\mathbf{A}^{-1}\|}{1-(1-\varepsilon)}=
\varepsilon^{-1}\bigl \|\mathbf{A}^{-1}\bigr\|.
\]
By using the last relation, we complete the proof,
\begin{eqnarray*}
\bigl\|\mathbf{B}^{-1}-\mathbf{A}^{-1}\bigr\| &\le &\|
\mathbf{A}^{-1}\| \bigl\|\mathbf{A}\mathbf{B}^{-1} -\mathbf{I}\bigr\|
\\
& \le &\bigl\|\mathbf{A}^{-1}\bigr\| \|\mathbf{A}-\mathbf{B}\| \bigl\|
\mathbf{B}^{-1}\bigr\| \le \varepsilon^{-1}\bigl\|\mathbf{A}^{-1}
\bigr\|^2 \|\mathbf{A}-\mathbf{B}\|.
\end{eqnarray*}
\upqed\end{pf}
\end{appendix}

\section*{Acknowledgments}
The comments of the referees and the
Associate Editor helped us to better organize the paper.


\printaddresses

\end{document}